\documentclass[graybox]{svmult}

\usepackage{mathptmx}       
\usepackage{helvet}         
\usepackage{courier}        
\usepackage{type1cm}        
\usepackage{makeidx}         
\usepackage{graphicx}        
\usepackage{multicol}        
\usepackage[bottom]{footmisc}

\usepackage{amssymb,amsmath}

\newcommand{\lavg}{{ \big\{\hspace{-0.99ex}\big\{ }}						
\newcommand{\ravg}{{ \big\}\hspace{-0.99ex}\big\} }}

\newcommand\avg[1]{{ \lavg#1\ravg }}

\newcommand{\tripnorm}[1]{{\left\vert\kern-\nulldelimiterspace\left\vert\kern-\nulldelimiterspace\left\vert #1
	\right\vert\kern-\nulldelimiterspace\right\vert\kern-\nulldelimiterspace\right\vert}}

\newcommand\restr[2]{{												
	\left.\kern-\nulldelimiterspace									
	#1
	\vphantom{\big|}
	\right|_{#2}
	}}

\makeindex             

\begin{document}

\title*{Implicit LES with high-order H(div)-conforming FEM for incompressible Navier-Stokes flows}
\titlerunning{Implicit LES with high-order H(div)-conforming FEM for incompressible flows} 
\author{Gert Lube \& Philipp W. Schroeder}
\institute{G.~Lube \at Georg-August University G\"ottingen, Institut for Numerical and Applied Mathematics,
           Lotzestrasse 16-18, D-37085 G\"ottingen, Germany, \email{lube@math.uni-goettingen.de}
\and P.W.~Schroeder \at Georg-August University G\"ottingen, Institut for Numerical and Applied Mathematics,
           Lotzestrasse 16-18, D-37085 G\"ottingen, Germany, \email{p.schroeder@math.uni-goettingen.de}}
%
%
\maketitle

\abstract{Consider the transient incompressible Navier-Stokes flow at high Reynolds numbers. 
A high-order H(div)-conforming FEM with pointwise divergence-free discrete velocities is applied to
implicit large-eddy-simulation in two limit cases: \\
i) decaying turbulence in periodic domains, 
ii) wall bounded channel flow. }

\section{H(div)-conforming dGFEM for Navier-Stokes problem}  \label{sec-1}
%
Consider a flow in a bounded polyhedron $\Omega \subset {\mathbb R}^d,~ d \le 3$ with boundary 
$\partial \Omega = \Gamma_0 \cup \Gamma_{per}$ and outer unit normal  ${\bf n}=(n_i)_{i=1}^d$. 
Set $Q_T := (0,T) \times \Omega$ and denote ${\bf f}$ as source term. 
We want to find velocity ${\bf u}: Q_T  \to {\mathbb R}^d$ and pressure $p: Q_T \to {\mathbb R}$ s.t.
\begin{eqnarray} \label{TINS1}
   \partial_t {\bf u}- \nu \Delta {\bf u} + ({\bf u} \cdot \nabla){\bf u} +  
        \nabla p & = & {\bf f} \qquad \mbox{in}~ Q_T , \\
  \label{TINS2}
       \nabla \cdot {\bf u} & = & 0 \qquad \mbox{in}~ Q_T , \\
   \label{TINS3}
      {\bf u}  & = & {\bf 0} \qquad \mbox{on}~ (0,T) \times \Gamma_0 , \\
   \label{TINS4}
     {\bf u} & = & {\bf u}_0 \quad \mbox{on}~ \lbrace 0 \rbrace \times  \Omega ,
\end{eqnarray}
and periodic boundary conditions on $\Gamma_{per}$. Let ${\bf H} = [L^2(\Omega)]^d$ with inner product  
$(\cdot,\cdot)_{\bf H}$ and assume ${\bf u}_0 \in {\bf H},~ {\bf f} \in L^2(0,T;{\bf H})$. 
The inner product in $L^2(\Omega)$ is $(\cdot, \cdot)_\Omega$. 

A variational formulation of the transient incompressible Navier-Stokes problem (\ref{TINS1})-(\ref{TINS4}) 
is to find $({\bf u},p) \in {\bf X} \times Q \subseteq [H^1(\Omega)]^d \times L^2(\Omega)$ for $t \in (0,T)$ 
{\it a.e.} from
\begin{eqnarray} \label{Nav1}
   (\partial_t {\bf u}(t),{\bf v})_{\bf H} + \nu a({\bf u}(t),{\bf v}) + c({\bf u}(t),{\bf u}(t),{\bf v})
    + b(p(t),{\bf v}) & = & ({\bf f}(t),{\bf v})_{\bf H},  \\  
   \label{Nav2}
                  -b(q,{\bf u}(t)) & = & 0,  \\
   \label{Nav3}
                        {\bf u}(0)   & = & {\bf u}_0 .
\end{eqnarray}
with bounded  bilinear resp. trilinear forms
\begin{equation} \label{forms}
   a({\bf u},{\bf v}):= \nu (\nabla {\bf u}, \nabla {\bf v})_\Omega,~~
           b(q,{\bf v}):= -(q, \nabla \cdot {\bf v})_\Omega,~~
     c({\bf w},{\bf u},{\bf v}):= (({\bf w} \cdot \nabla) {\bf u},{\bf w})_\Omega .
\end{equation}

Consider ${\bf H}(\text{div})$-conforming, discontinuous Galerkin methods (dGFEM) 
with
\begin{eqnarray} \label{Hdiv}
      {\bf H}(\text{div};\Omega) & := & \lbrace {\bf w} \in {\bf H}:~~ 
                  \nabla \cdot {\bf w} \in L^2(\Omega) \rbrace, \\
    \label{Hdiv0}
  {\bf H}_{\Gamma_0}(\text{div};\Omega) &:= &\lbrace {\bf v} \in {\bf H}(\text{div};\Omega):~ 
              {\bf v} \cdot {\bf n}|_{\Gamma_0} = 0 \rbrace .
\end{eqnarray}
Let ${\mathcal T}_h$ be a shape-regular decomposition of $\Omega \subset {\mathbb R}^d$. Moreover,
denote ${\mathcal E}_h$ the set of (open) edges ($d=2$) or faces ($d=3$) in ${\mathcal T}_h$.
${\mathcal E}^B_h \subset {\mathcal E}_h$ is the set of all $E \in {\mathcal E}_h$ with 
$E \cap \Gamma_0 \ne \emptyset$ and ${\mathcal E}_h^I :={\mathcal E}_h \setminus {\mathcal E}^B_h$
the set of interior edges. Please note that edges/faces on $\Gamma_{per}$ are considered as interior 
edges/faces. 
Consider adjacent elements $K, K' \in {\mathcal T}_h$ with $\partial K \cap \partial K' = E$ and 
unit normal vector ${\bf \mu}_E$. 
For a scalar function $v$ in the broken Sobolev space $H^{1}(\Omega,\mathcal{T}_h)$ denote jump resp. 
average of $v$ across $E$ by 
\begin{equation}
   [|v|]_E := v|_{\partial K\cap E}-v|_{\partial K'\cap E}, \qquad 
  \avg{v}_E : = (v|_{\partial K\cap E}+v|_{\partial K'\cap E})/2. 
\end{equation}
For ${\bf v} \in [H^{1}(\Omega,\mathcal{T}_h)]^d$, jump and average are understood component-wise.
\begin{lemma}\cite{diPietro-Ern}~
    Let ${\bf W}_h$ be a space of vector-valued polynomials w.r.t. ${\mathcal T}_h$.
    Then ${\bf W}_h \subset {\bf H}(\text{div};\Omega)$ if $[|{\bf v}|]_E \cdot {\bf \mu}_E= 0$ 
    for all ${\bf v} \in {\bf W}_h$ and all $E \in {\mathcal E}_h^I$.
\end{lemma}
Owing to Lemma 1 $[|{\bf v}|]_{t,E}= [|{\bf v}-({\bf v}\cdot \mu_E)\mu_E|]_E$ is the tangential jump 
across $E \in {\mathcal E}_h$.
\begin{example}~Examples of ${\bf H}(\text{div})$-conforming FEM are given in \cite{Boffi-etal}.
On simplicial grids one can apply Raviart-Thomas (RT) or Brezzi-Douglas-Marini (BDM) spaces
 \begin{eqnarray} \label{15.3a}
   \mbox{RT}_k & = &\lbrace {\bf w}_h \in {\bf H}_{\Gamma_0}(\mbox{div};\Omega):~ 
    {\bf w}_h|_K \in {\mathbb P}_k(K) \oplus x {\mathbb P}_k(K)~  \forall K \in {\mathcal T}_h \rbrace,~
    k \in {\mathbb N}_0~ \\
    \label{15.4}
     \mbox{BDM}_k & = & \lbrace {\bf w}_h \in {\bf H}_{\Gamma_0}(\text{div};\Omega):~ 
         {\bf w}_h|_K \in {\mathbb P}_k(K)~ \forall K \in {\mathcal T}_h \rbrace,~ k \in {\mathbb N}. 
  \end{eqnarray}
On quadrilateral meshes, local Raviart-Thomas (RT) elements of degree $k \in {\mathbb N}_0$ are
$\mbox{RT}_k(K) = ({\mathbb P}_{k+1,k}(K),{\mathbb P}_{k,k+1}(K))^t, d=2$. For $d=3$, one has similarly
$\mbox{RT}_k(K)  = ({\mathbb P}_{k+1,k,k}(K),{\mathbb P}_{k,k+1,k}(K),{\mathbb P}_{k,k,k+1}(K))^t$.
        \hfill $ \Box $
\end{example}
Let ${\bf w}_h \in {\bf W}_h \subset {\bf H}(\mbox{div};\Omega)$ with 
${\bf W}_h  \in \lbrace \mbox{RT}_k, \mbox{BDM}_k \rbrace$.
The spaces ${\bf W}_h \not\subset [H^1(\Omega)]^d$ are not $[H^1(\Omega)]^d$-stable, hence not 
directly applicable to the Navier-Stokes problem. 
As a remedy, we modify the diffusion bilinear form $a$ using a symmetric interior penalty (SIP) 
dGFEM-approach with the broken gradient $\nabla _h {\bf v} := \nabla ({\bf v}|_K)$:
For sufficiently smooth ${\bf u} \in [H^s(\Omega)]^d, s>\frac32$, we define by adding two consistent 
terms 
\begin{eqnarray} \label{a-mod} 
  a_h({\bf u},{\bf w}_h) & := & \int_\Omega \nabla_h {\bf u}:\nabla_h {\bf w}_h~dx 
     + \sum_{E \in {\mathcal E}_h} \sigma h_E^{-1} \int_E |[{\bf u}|]_t |[{\bf w}]|_t~ ds \\
  & -  & \nonumber
    \sum_{E \in {\mathcal E}_h} \int_{E} \big( \avg{\nabla_h {\bf u} \cdot \mu_E} [|{\bf w}_h|]_t
         + \avg{\nabla_h {\bf w}_h \cdot \mu_E} [|{\bf u}|]_t~\big)~ ds~~  \forall {\bf w}_h \in {\bf W}_h
\end{eqnarray} 
with $h_E := \mbox{diam}(E)$ and parameter $\sigma >0$ (to be chosen according to next lemma).
\smallskip \\
Define the following  discrete $H^1$-norms  $\|{\bf w}\|_{1,h}$ and  $\|{\bf w}\|_{1,h,*}$ 
\begin{eqnarray} \label{H1-norm1}
     \|{\bf w}\|^2_{1,h} & := & \sum_{K \in {\mathcal T}_h} \| \nabla {\bf w}\|^2_{L^2(K)} 
          + \sum_{E \in {\mathcal E}_h} h_E^{-1} \| [|{\bf w}|]_\tau\|^2_{L^2(E)} , \\
    \label{H1-norm2}
    \|{\bf w}\|^2_{1,h,*} & := & \|{\bf w}\|^2_{1,h} +  \sum_{E\in {\mathcal E}_h}
          h_E \| \avg{\nabla _h {\bf w} \cdot \mu_E}\|^2_{L^2(E)}.
    \end{eqnarray}    
\begin{lemma} \cite{diPietro-Ern}
   There exists constant $\sigma_0$ (depending only on $k$ and on shape regularity of 
    ${\mathcal T}_h$) s.t. for $\sigma \ge \sigma_0$ one has:
    \begin{eqnarray} \label{15.11} 
       a_h({\bf w}_h,{\bf w}_h) & \ge & \frac12 \|{\bf w}\|_{1,h}^2 \quad \forall {\bf w}_h \in {\bf W}_h, \\
     \label{15.13} 
       a_h({\bf v},{\bf w}_h) & \le & C \|{\bf v}\|_{1,h,*} \|{\bf w}_h\|_{1,h} \quad
       \forall {\bf w}_h \in {\bf W}_h ~~ \mbox{and}~~ {\bf v} \in [H^s(\Omega)]^d,~ s > \frac32 .
    \end{eqnarray}
\end{lemma}
\begin{lemma} \cite{Boffi-etal}
\text{RT}- and \text{BDM}-spaces, together with appropriate discrete spaces $Q_h$
\begin{eqnarray*}
  {\bf W}_h & = & \mbox{RT}_k ~~\mbox{with}~~
     Q_h := \lbrace q_h \in L^2(\Omega):~ q_h|_K \in {\mathbb P}_k(K)~ \forall K \in {\mathcal T}_h \rbrace 
  \qquad \text{and}  \\
 {\bf W}_h & = & \mbox{BDM}_k ~~\mbox{with}~~
      Q_h := \lbrace q_h \in L^2(\Omega):~ q_h|_K \in {\mathbb P}_{k-1}(K)~ \forall K \in {\mathcal T}_h \rbrace
\end{eqnarray*}
form inf-sup stable pairs w.r.t. discrete $H^1$-norm:
\begin{equation} \label{15.16}
       \exists \beta_h \ge \beta_0 > 0~ \mbox{s.t.}~~ 
       \inf_{q_h \in Q_h \setminus \lbrace 0 \rbrace } \sup_{{\bf w}_h \in {\bf W}_h \setminus \lbrace {\bf 0} \rbrace }
          \frac{ (\nabla \cdot {\bf w}_h,q_h)_\Omega}{\| {\bf w}_h\|_{1,h} \| q_h\|_{L^2(\Omega)}} \ge \beta_h .
\end{equation}
By construction~ $\nabla \cdot {\bf W}_h = Q_h$, these spaces are globally pointwise 
divergence-free:
\begin{equation} \label{divfree}
    \lbrace {\bf w}_h \in {\bf W}_h:~ (\nabla \cdot {\bf w}_h,q_h)_\Omega =0~ \forall q_h \in Q_h \rbrace
    = \lbrace {\bf w}_h \in {\bf W}_h:~ \nabla \cdot {\bf w}_h = 0 \rbrace . 
\end{equation}
\end{lemma}
For an exactly divergence-free field ${\bf b} \in [L^\infty(\Omega)]^d \cap {\bf H}(\text{div};\Omega)$
we modify the convective term $c$ as in \cite{diPietro-Ern}  by
  \begin{eqnarray} \nonumber
    c_h({\bf b};{\bf u},{\bf v}) & := & \sum_{K \in {\mathcal T}_h} (({\bf b} \cdot \nabla ){\bf u},{\bf v})_K \\
    \label{convec} 
     & - &  \sum_{E \in {\mathcal E}_h^i} \big( ({\bf b} \cdot {\bf \mu}_E) ([|{\bf u}|], \avg{\bf v}) \big)_E
        + \frac12 \sum_{E \in {\mathcal E}_h^i} \big( |{\bf b} \cdot {\bf \mu}_E| [|{\bf u}|],[|{\bf v}|]) \big)_E.
  \end{eqnarray}
The first right-hand side terms corresponds to the standard form of the convective term. 
The last two facet terms, the {\it upwind} discretization, are consistent perturbations of the 
standard form of the convective term for ${\bf u},{\bf v} \in {\bf X}$.
The impact of these terms is included in the jump semi-norm  $|{\bf v}|_{{\bf b},\text{upw}}$
defined via 
\begin{equation} \label{upwind}
      |{\bf v}|^2_{{\bf b},\text{upw}} := \frac12 \sum_{E \in {\mathcal E}_h^i} 
                 |{\bf b} \cdot {\bf \mu}_E|~ \| [|{\bf v}|] \|_{L^2(E)}^2.
\end{equation} 
In case of exactly divergence-free fields ${\bf b}$, one has~
     $c_h({\bf b};{\bf v},{\bf v}) =  |{\bf v}|^2_{{\bf b},\text{upw}}.$

We consider now the ${\bf H}(\text{div})$-conforming dGFEM for the transient Navier-Stokes
problem (\ref{Nav1})-(\ref{Nav3}) with  ${\bf f} \in L^2(0,T;{\bf H})$. Combining 
the SIP-form of the diffusive term and the upwind-discretization of the convective term, one obtains: \\
Find $({\bf u}_h,p_h): (0,T) \to {\bf W}_h \times Q_h$ with ${\bf u}_h(0)={\bf u}_{0,h}$ s.t.
for all $({\bf v}_h,q_h) \in {\bf W}_h \times Q_h$:
 \begin{eqnarray}  \label{Hdiv-FEMa}
  (\partial_t {\bf u}_h,{\bf v}_h)_{\bf H} + \nu a_h({\bf u}_h,{\bf v}_h) +
   c_h({\bf u}_h;{\bf u}_h,{\bf v}_h)    
          +   b(p_h,{\bf v}_h) & = &  ({\bf f},{\bf v}_h)_{\bf H} , \\
 \label{Hdiv-FEMb}
     - b(q_h,{\bf u}_h) & = & 0 . 
 \end{eqnarray}

All computations have been done using a hybridized variant of (\ref{Hdiv-FEMa})-(\ref{Hdiv-FEMb})
implemented in the high-order software package NGSolve \cite{Schoeberl}.

We will consider method (\ref{Hdiv-FEMa})-(\ref{Hdiv-FEMb}) as tool for implicit large-eddy-simulation 
(ILES) in two limit cases: 
i) decaying turbulence in periodic 2D and 3D domains (see Sec.~\ref{sec-2}) and
ii) wall bounded flow in a 3D-channel (see Sec.~\ref{sec-3}).

\section{Decaying 2D- and 3D-turbulent flows} \label{sec-2}
%

\subsection{Stability and error analysis for decaying flows} \label{subsec:2D}
%
Consider now decaying flows, i.e. we consider problem (\ref{Hdiv-FEMa})-(\ref{Hdiv-FEMb}) with 
${\bf f} \equiv {\bf 0}$. 
Using the mesh-dependent expressions (\ref{H1-norm1}) and (\ref{upwind})
and setting ${\bf v}_h={\bf u}_h$ in the semidiscrete problem (\ref{Hdiv-FEMa})-(\ref{Hdiv-FEMb}), one 
obtains with $\| {\bf v} \|_e^2 := a_h({\bf v},{\bf v})$ the balance 
\begin{equation} \label{balance}
 \frac{d}{dt} \left( \frac12 \|{\bf u}_h(t)\|_{L^2(\Omega)}^2 \right) + \nu \| {\bf u}_h(t)\|^2_e + 
             |{\bf u}_h(t)|_{{\bf u}_h,\text{upw}}^2    =   0  .
\end{equation}
This implies existence of 
$({\bf u}_h,p_h)$ and bounds for kinetic and dissipation energies:
\begin{eqnarray} \label{stab1}
 \frac12 \| {\bf u}_h(t) \|_{L^2(\Omega)}^2 & \le &
    \frac12 \|{\bf u}_{0h}\|_{L^2(\Omega)}^2  \exp (-\nu t/C_F^2) , \\
 \label{stab2}
 \int_0^t \Big( \frac{\nu}{2} \|{\bf u}_h(\tau)\|_e^2 + |{\bf u}_h(\tau)|_{{\bf u}_h,\text{upw}}^2 \Big) d\tau 
        & \le & \frac12 \|{\bf u}_{0h}\|_{L^2(\Omega)}^2  .
\end{eqnarray}
In case of smooth velocity with ${\bf u} \in L^1(0,T;[W^{1,\infty}(\Omega)]^d)$, we obtain the following 
pressure-robust and $Re$-semi-robust error estimate.
\begin{theorem} \cite{Schroeder-18}   \label{theorem:converg}
Let ${\bf u} \in L^2(0,T;{\bf H}^{\frac32+\epsilon}(\Omega)),~ \epsilon>0$, 
$\nabla {\bf u} \in L^1(0,T,[L^\infty(\Omega)]^d)$ and ${\bf u}_h(0)= \pi_S{\bf u}_0$ with 
Stokes projector $\pi_S {\bf u}$, i.e. $a_h({\bf u}-\pi_S{\bf u},{\bf v}_h)=0~ \forall 
{\bf v}_h \in {\bf W}_h$; then:
\begin{eqnarray*}
 & & \frac12 \| {\bf u}_h- \pi_S {\bf u}\|^2_{L^\infty(0,T;L^2(\Omega))} 
   + \int_0^T \big[ \frac{\nu}{2} \|{\bf u}_h- \pi_S {\bf u}\|_{1,h}^2 + 
    | {\bf u}_h-\pi_S {\bf }|_{{\bf u}_h,\mbox{\scriptsize upw}}^2 \big] ~d\tau   \\
 & \le & e^{G_{\bf u}(T)}  \int_0^T \big[ \| \partial_t \eta \|_{L^2(\Omega)}^2 +  
     \| {\bf u}\|_{L^\infty(\Omega)} \| \nabla _h \eta \|^2_{L^2(\Omega)}  
         + h^{-2} \| \nabla {\bf u} \|_{L^\infty(\Omega)} \| \eta \|^2_{L^2(\Omega)} \big] ~d\tau
\end{eqnarray*}
with $\eta := {\bf u}-\pi_S{\bf u}$ and Gronwall factor
$$G_{\bf u}(T)  := T + \| {\bf u} \|_{L^1(0,T;[L^\infty (\Omega)]^d)} + 
               C \| \nabla {\bf u} \|_{L^1(0,T;[L^\infty (\Omega)]^d)} .$$
\end{theorem}

The vorticity equation for $\omega := \nabla \times {\bf u}$ describes the dynamics of 
decaying flows:
\begin{equation} \label{vorticity}
  \partial_t \omega + {\bf u} \cdot \nabla \omega - \nu \Delta \omega = 
       \omega \cdot \nabla {\bf u}, \qquad \nabla \cdot {\bf u} = 0.
\end{equation}
The vortex stretching term~ $\omega \cdot \nabla {\bf u}$ vanishes for $d=2$ which leads to a 
completely different behavior for $d=2$ and $d=3$.

\subsection{Decaying 2D-turbulent flow} \label{subsec:2D}
%
Consider the following 2D-turbulent flow problem with a unique solution of (\ref{Nav1})-(\ref{Nav3}).
\begin{example}~~ 2D-lattice flow \\
Consider on $\Omega=(-1,1)^2$ the following solution of the steady Euler model ($\nu=0$)
\[
  {\bf u}_0(x)  = (- \Psi_{x_2}(x),\Psi_{x_1}(x))^t, \qquad
  \Psi (x) := \frac{1}{2\pi} \sin (2\pi x_1)\cos (2\pi x_2) .
\]
The initial vorticity $\omega_0 = \nabla \times {\bf u}_0$ is shown in Fig.~\ref{fig:ICLattice-structure}
for $t=0$. The Taylor cells ${\bf u}(t,x) = {\bf u}_0(x) e^{-4\pi^2 \nu t}$ are the (unique!) solution 
of the transient Navier-Stokes model. For this very smooth solution, a high-order FEM is preferable. 

This is a generalized Beltrami flow, since~ 
$({\bf u} \cdot \nabla) {\bf u} = -\nabla p$. 
Due to pressure-robustness, a linearization via dropping $({\bf u} \cdot \nabla ) {\bf u}$ 
preserves the {\it coherent} structures of the initial solution \cite{Gauger-2018}.
\begin{figure}[h]
\centering
\includegraphics[width=0.245\textwidth]{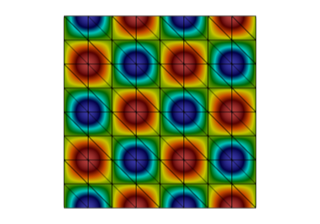} 
\includegraphics[width=0.245\textwidth]{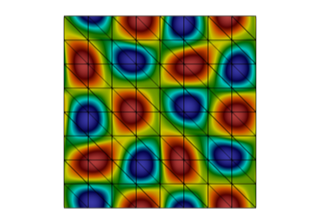} 
\includegraphics[width=0.245\textwidth]{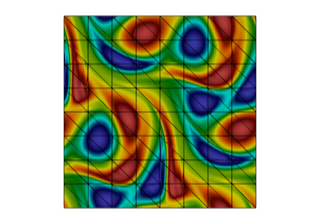} 
\includegraphics[width=0.245\textwidth]{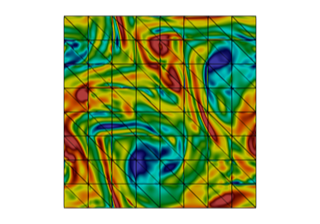} 
\\
\centering
\includegraphics[width=0.245\textwidth]{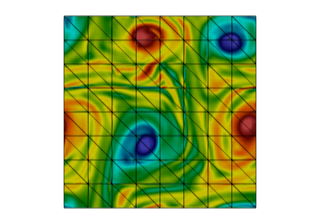} 
\includegraphics[width=0.245\textwidth]{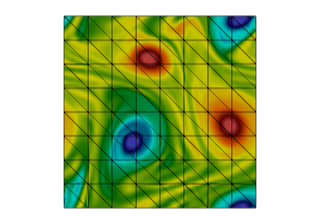} 
\includegraphics[width=0.245\textwidth]{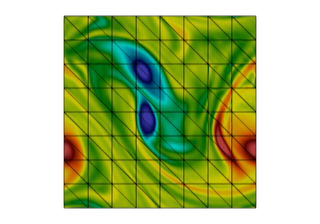} 
\includegraphics[width=0.245\textwidth]{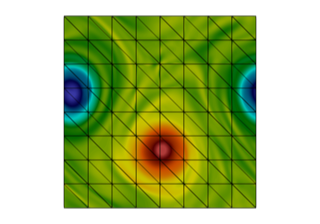} 
\caption{Example 2: Snapshots of vorticity $\omega_h = \nabla \times {\bf u}_h$~ of 
         high-order FEM with $k=8, N=8$ with  $t \in \lbrace 0, 22, 23, 26 \rbrace$ (see
         first row) and $t \in \lbrace 30, 35, 40, 50 \rbrace$ (see second row)}
\label{fig:ICLattice-structure}
\end{figure}
For order $k=8$ and $h=\frac14$, Fig.~~\ref{fig:ICLattice-structure} shows snapshots of the 
discrete vorticity on the time interval $0 \le t \le T=50$ for $\nu = 10^{-6}$. We observe 
a self-organization of vortical structures which deviates from the unique solution.

\begin{figure}[h]
\centering
\includegraphics[width=0.32\textwidth]{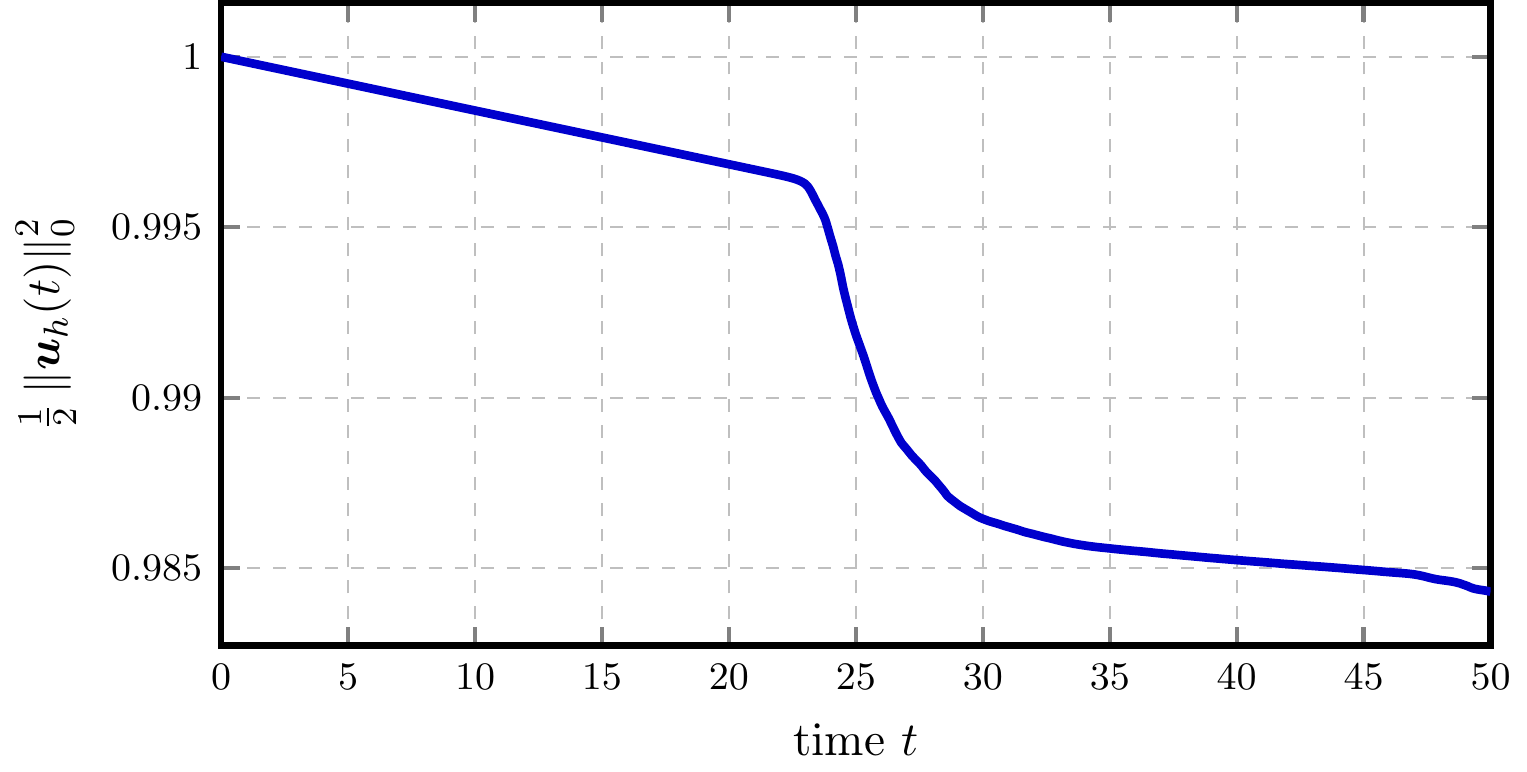} 
\includegraphics[width=0.32\textwidth]{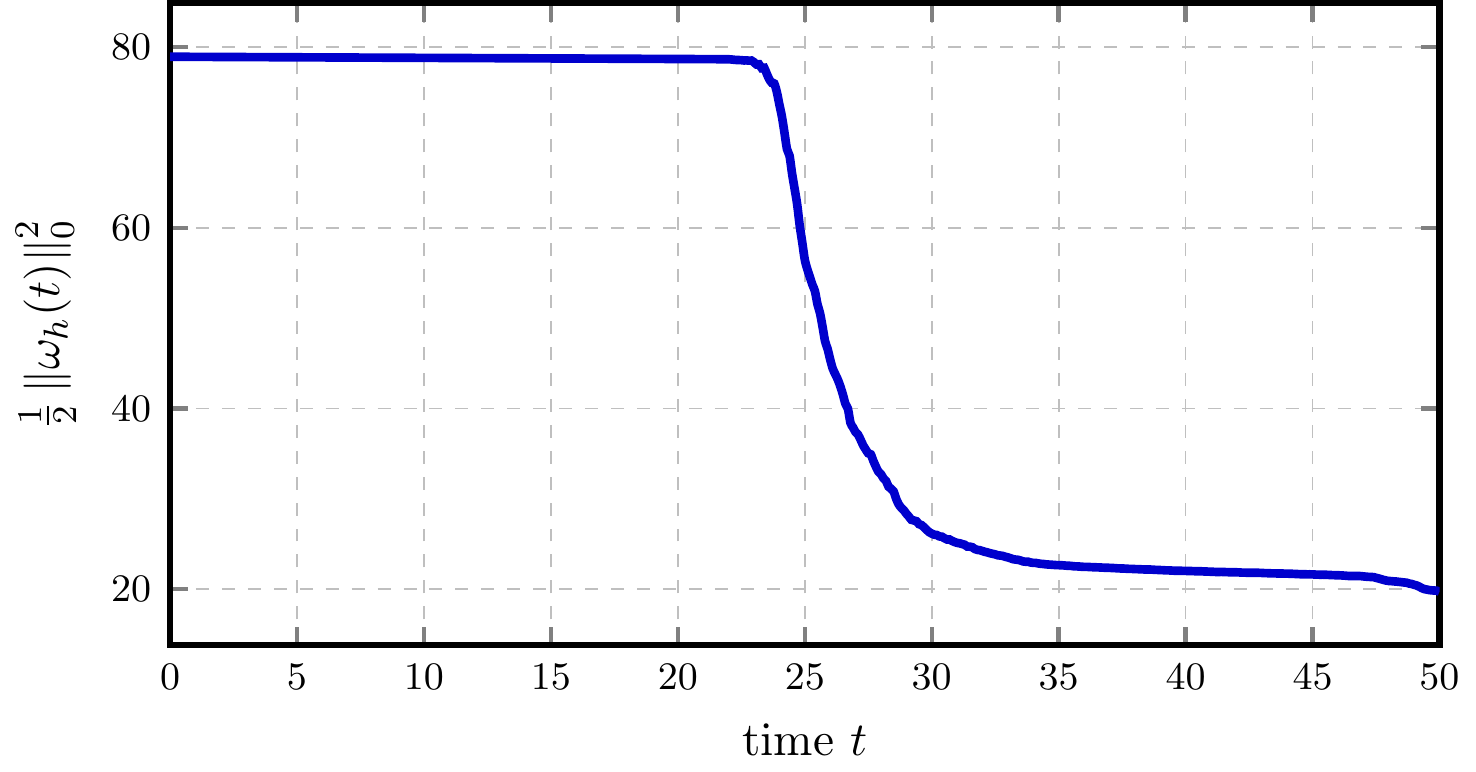} 
\includegraphics[width=0.32\textwidth]{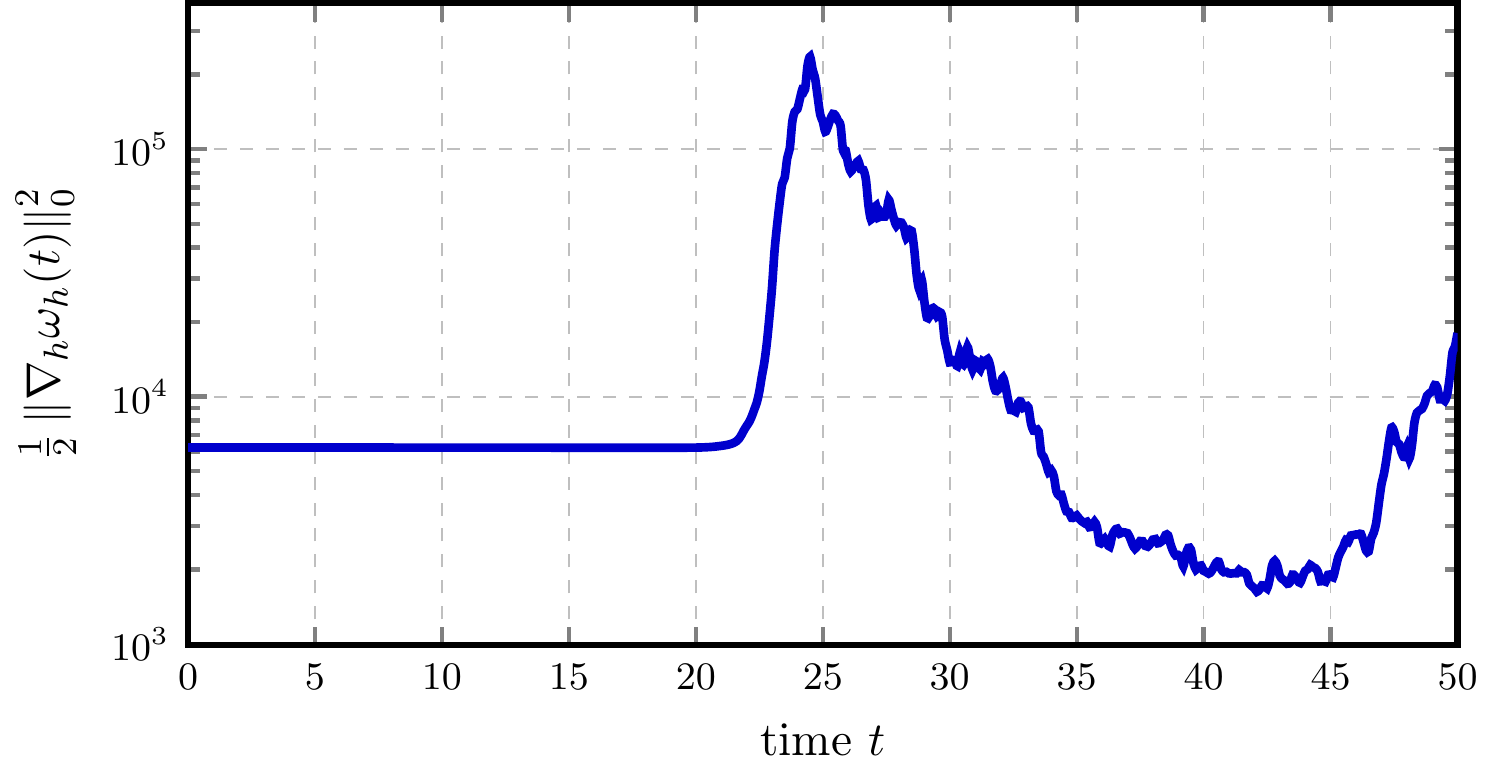} \\
\caption{Examle 2: Temporal development of kinetic energy, enstrophy and palinstrophy}
\label{fig:ICLatticesnap}
\end{figure}
Consider now the behavior of the kinetic energy $\frac12\|{\bf u}_h\|^2_{L^2(\Omega)}$, 
enstrophy $\frac12 \|\omega_h\|^2_{L^2(\Omega)}$ and palinstrophy $\frac12 \|\nabla_h \omega_h\|^2_{L^2(\Omega)}$
for $0\le t \le 50$, see Fig.~\ref{fig:ICLatticesnap}. Around $t=22.0$ the solution deviates from 
coherent structures of the exact solution, also visible in the strong reduction of the amplitude 
of the kinetic energy. 
\begin{figure}[h]
\centering
\includegraphics[width=0.35\textwidth]{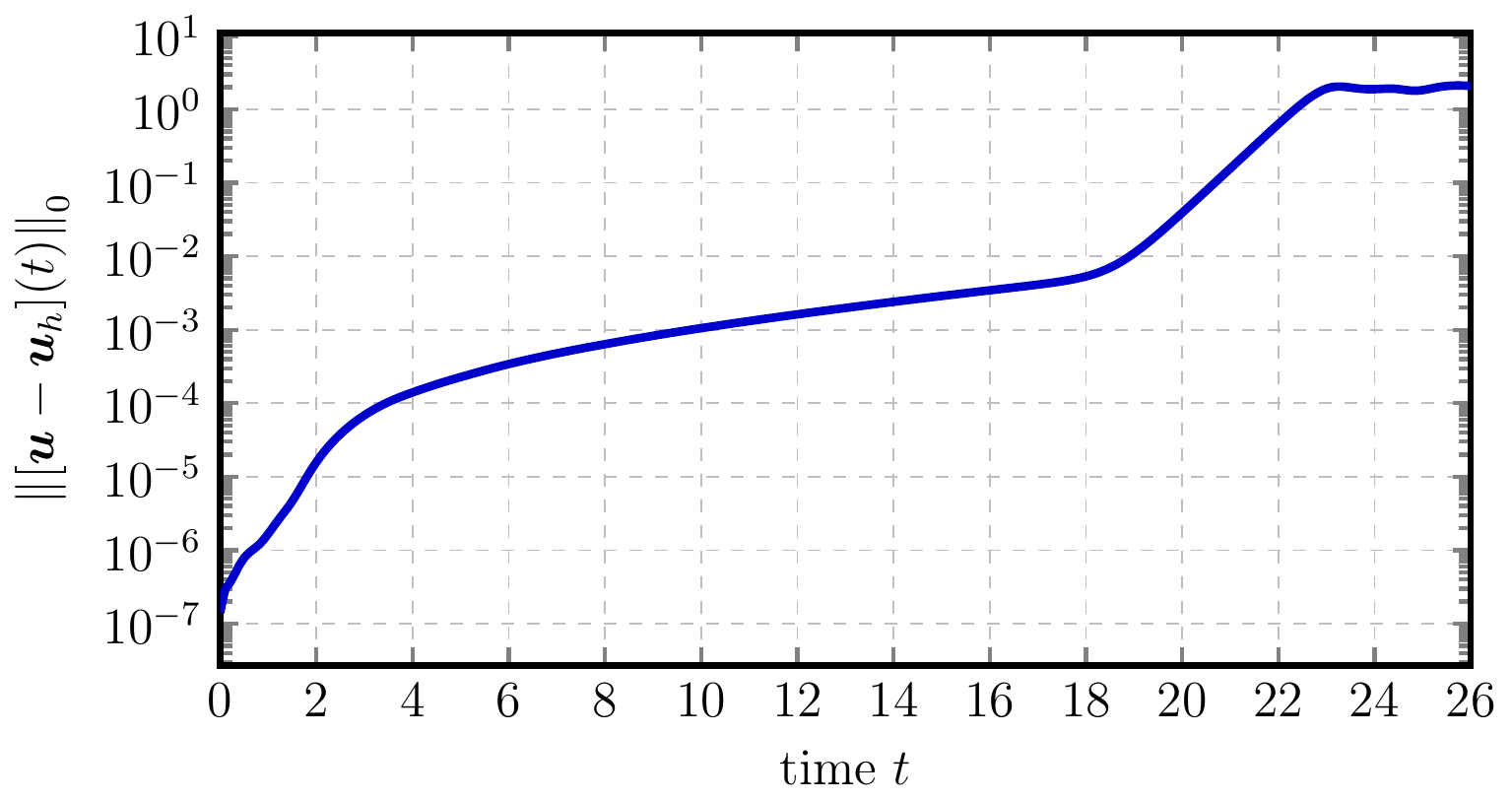}
\quad
\includegraphics[width=0.35\textwidth]{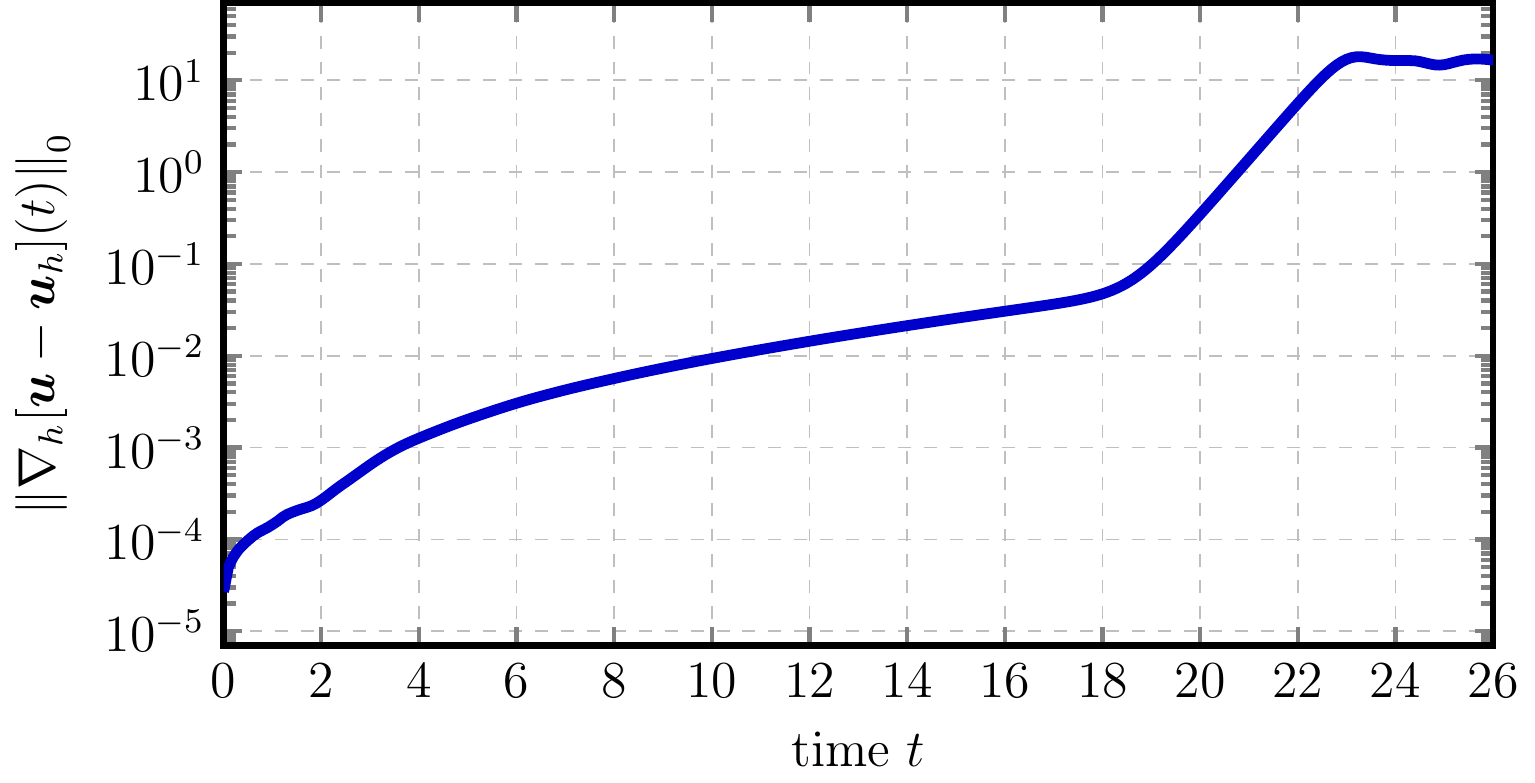} \\ 
\caption{Example 2: Error plots of high-order FEM for $\nu=10^{-6},~ k=8,~ N=8,~ \Delta t =10^{-3}$}
\label{fig:ICLatticeErrors}
\end{figure}
The exponential growth of the $L^2$- and $H^1$-errors of the velocity (according to 
Thm.~\ref{theorem:converg}) is shown in Fig.~\ref{fig:ICLatticeErrors}.
The initial condition ${\bf u}_0$ of the planar lattice flow induces a flow structure which, due 
to its saddle point structure, is {\it ''dynamically unstable so that small perturbations result 
in a very chaotic motion''} as stated in Majda \& Bertozzi \cite{Majda-2002}. 
A convincing discussion of self-organization in 2D-flows is given by van Groesen \cite{vanGroesen1988}.  

Note that the preservation of the coherent structures (of the unique solution) can be 
extended in time by higher order $k$ and/or $h$-refinement. Moreover, compared to standard 
mixed non-pressure-robust FEM, the application of pressure-robust FEM leads to much longer existence 
of such structures, see \cite{Gauger-2018}. \hfill $ \Box $
\end{example}
\begin{remark}
(i) A similar behavior of 2D-decaying turbulent flows is known for the 2D Kelvin-Helmholtz instability.
We refer to careful numerical studies in \cite{KH-2018}. 

(ii) The smallest scales depend on $d$. For $d=3$, one has Kolmogorov-length 
$\lambda_{3D} \approx L Re^{-\frac34}$ whereas for $d=2$, the Kraichnan-length is 
$\lambda_{2D} \approx L Re^{-\frac12}$. 
As conclusion, a direct numerical simulation (DNS) of 2D-flows at $Re \gg 1$ is much more 
realistic than in 3D, see \cite{John-2016}.  \hfill $ \Box $
\end{remark}

\subsection{Decaying 3D-turbulent flows} \label{subsec:3D}
%
From the vorticity equation (\ref{vorticity}) we concluded a completely different behavior of
high Re-number flows for $d=3$ as compared to $d=2$. The following example highlights the effect
of vortex stretching term $(\omega \cdot \nabla ){\bf u}$.
\begin{example}~ 3D-lattice flow \\ 
Consider the exact solution of the transient incompressible Navier-Stokes problem 
\[
      {\bf u}(t,x) = {\bf u}_0(x)  e^{-4\pi^2 \nu t}, \quad 
  {\bf u}_0(x_1,x_2) = (-\Psi_{x_2},\Psi_{x_1},\sqrt{2} \Psi )^t(x_1,x_2) 
\]
in $\Omega = (0,1)^3$ with stream function $\Psi$ as in Example~2, with ${\bf f}={\bf 0}$ and 
$\frac{1}{\nu}=2,000$. This problem can be seen as 3D-extension of the 2D-lattice flow 
\cite{Majda-2002}.

\begin{figure}[ht]
\begin{center}
\includegraphics[width=0.3\textwidth]{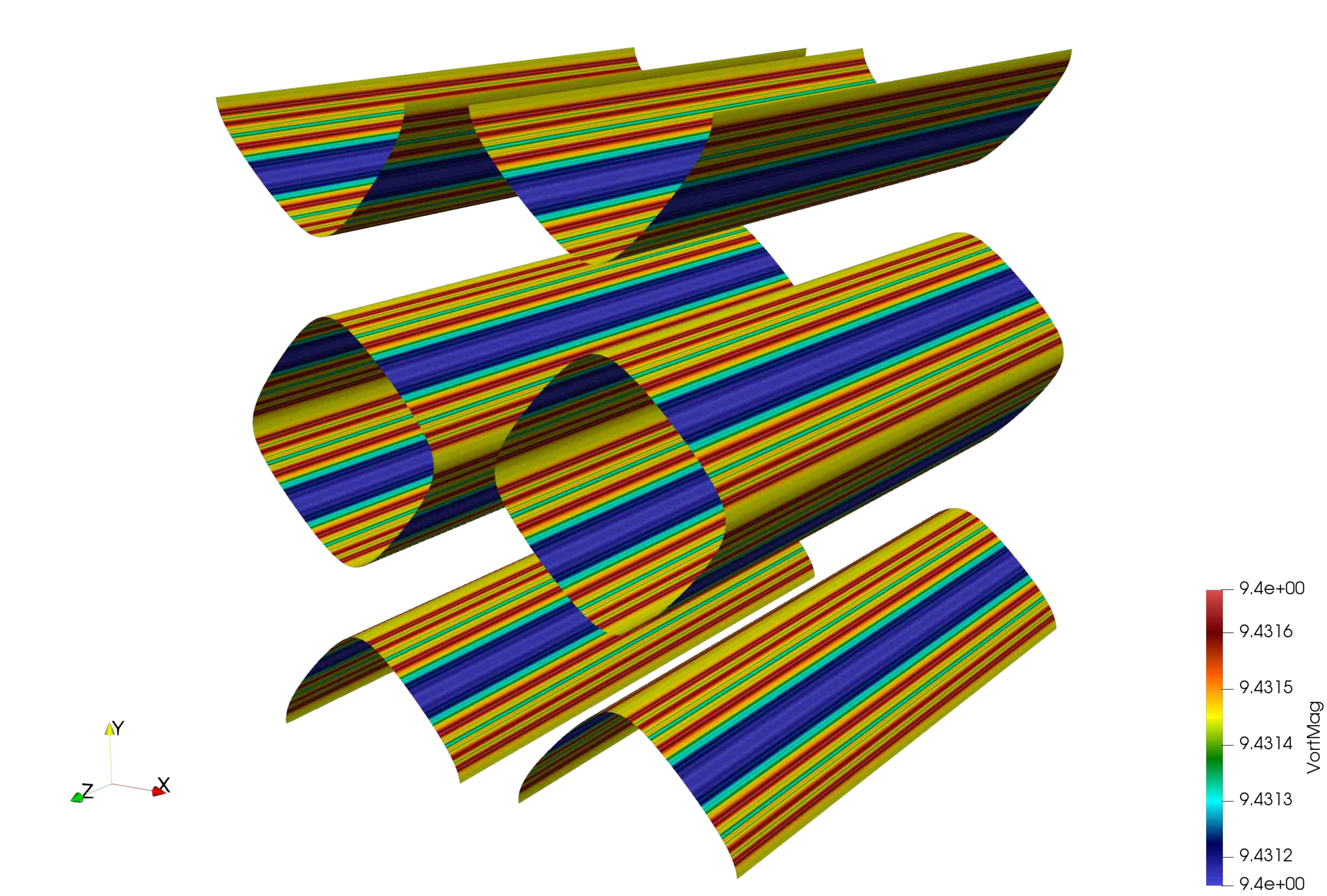} \quad 
\includegraphics[width=0.3\textwidth]{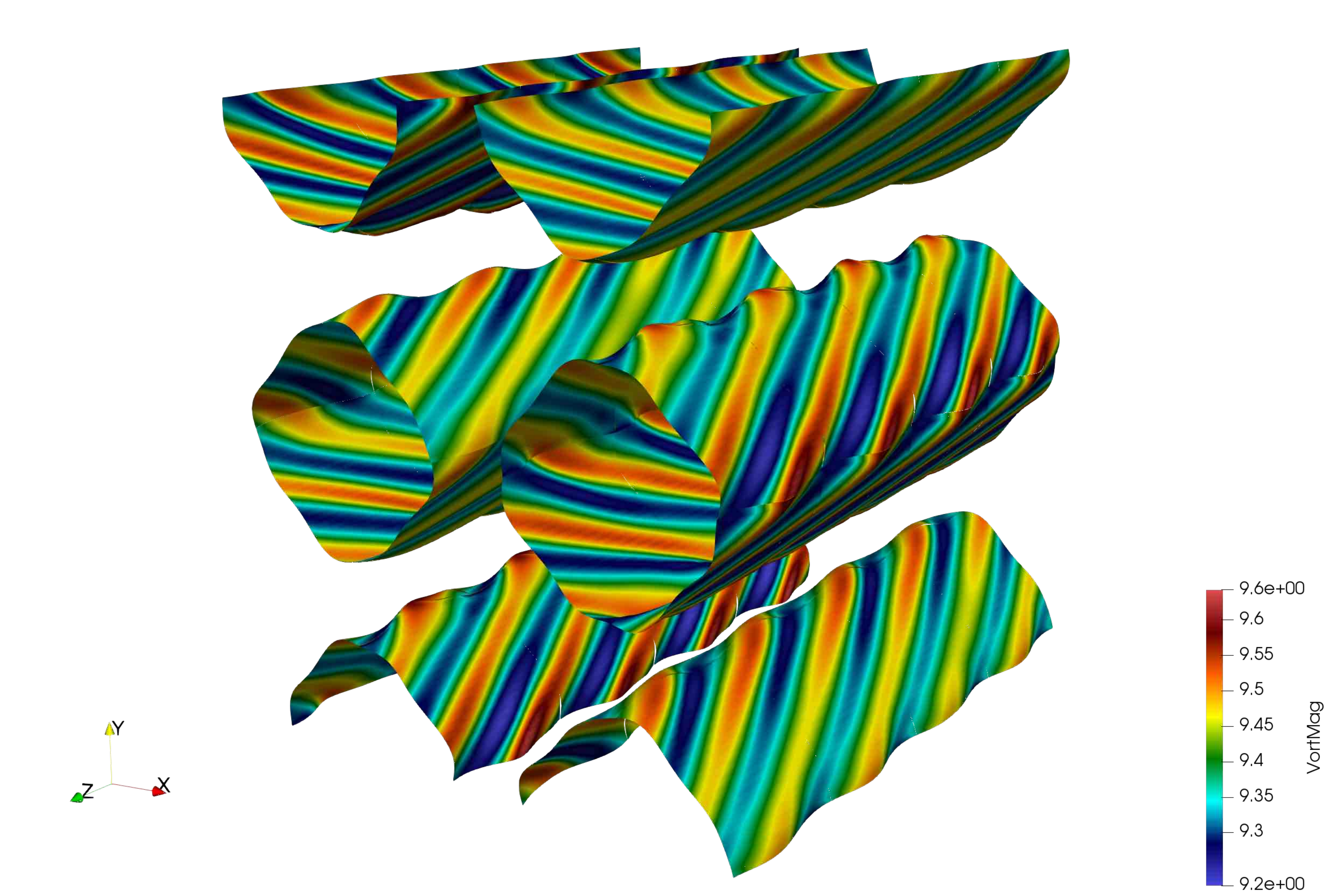} \quad 
\includegraphics[width=0.3\textwidth]{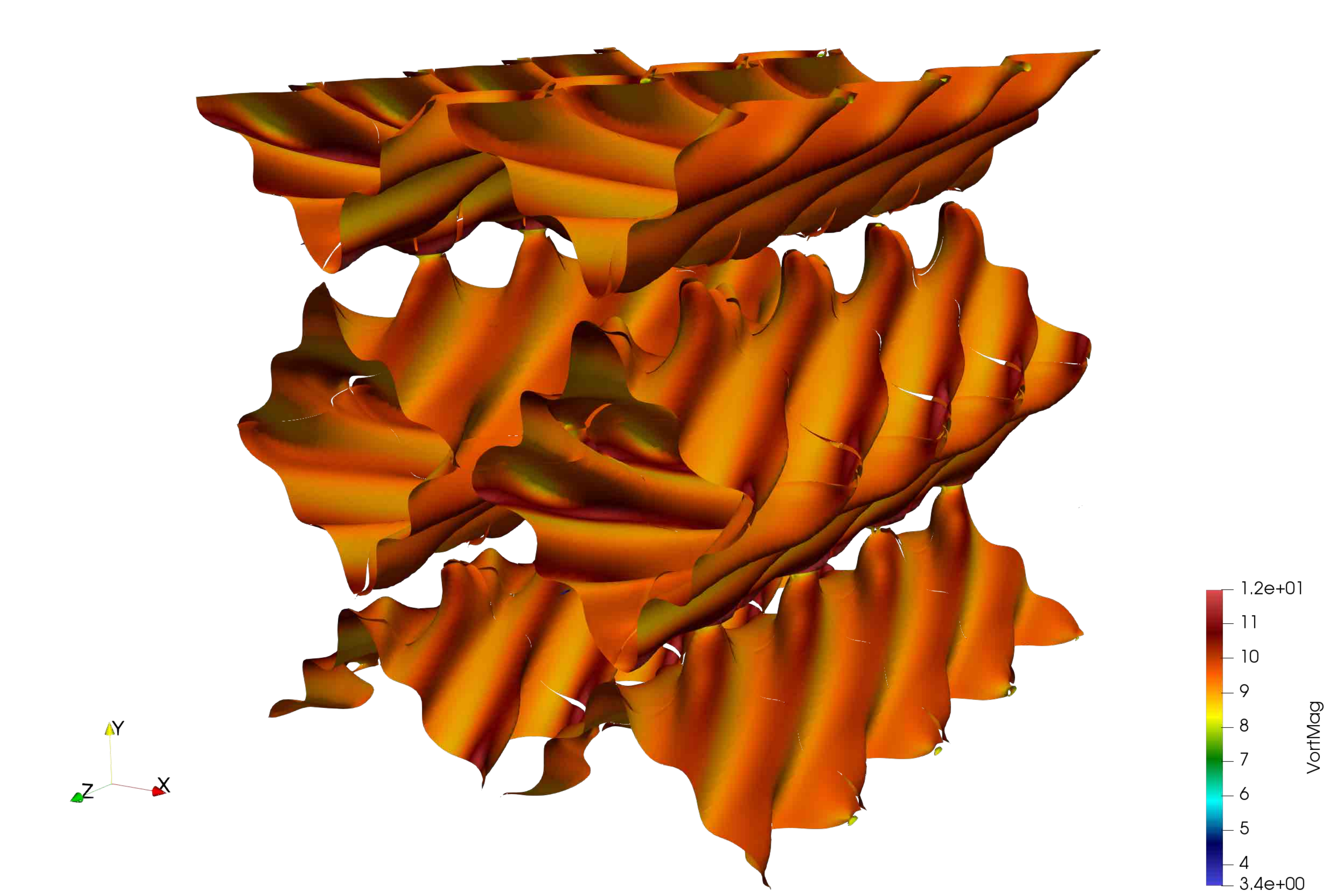}  \\
\includegraphics[width=0.3\textwidth]{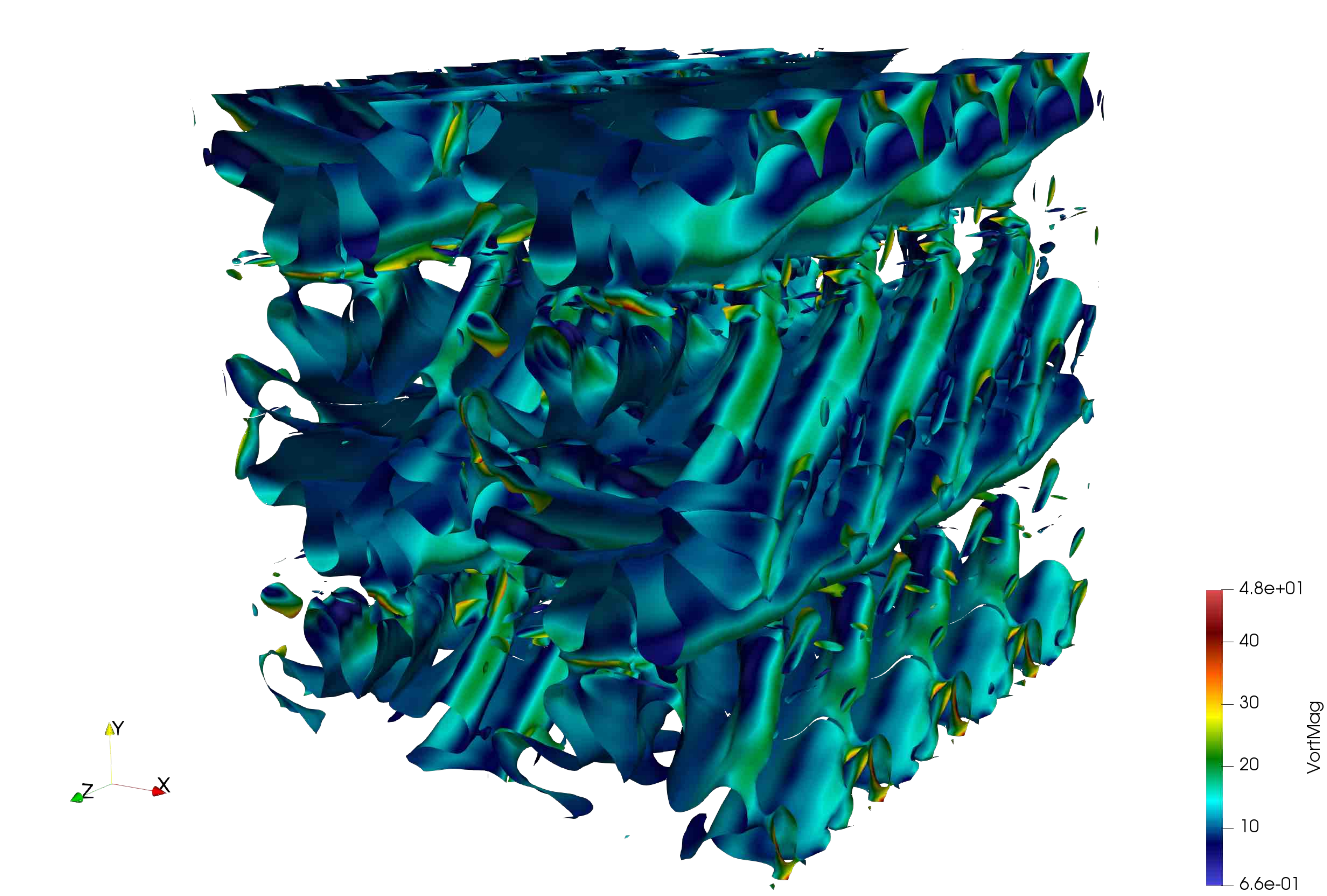} \quad
\includegraphics[width=0.3\textwidth]{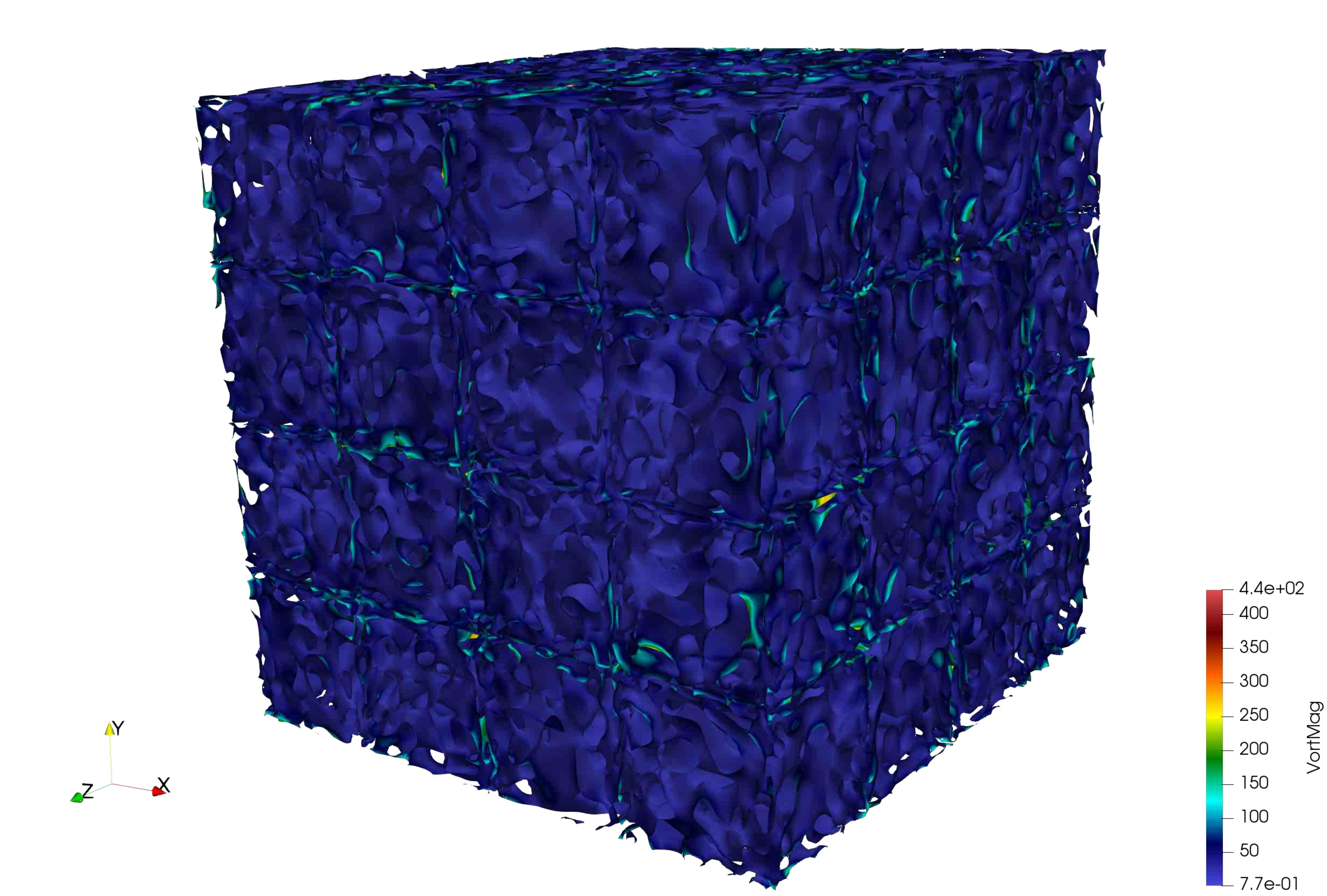} \quad
\includegraphics[width=0.3\textwidth]{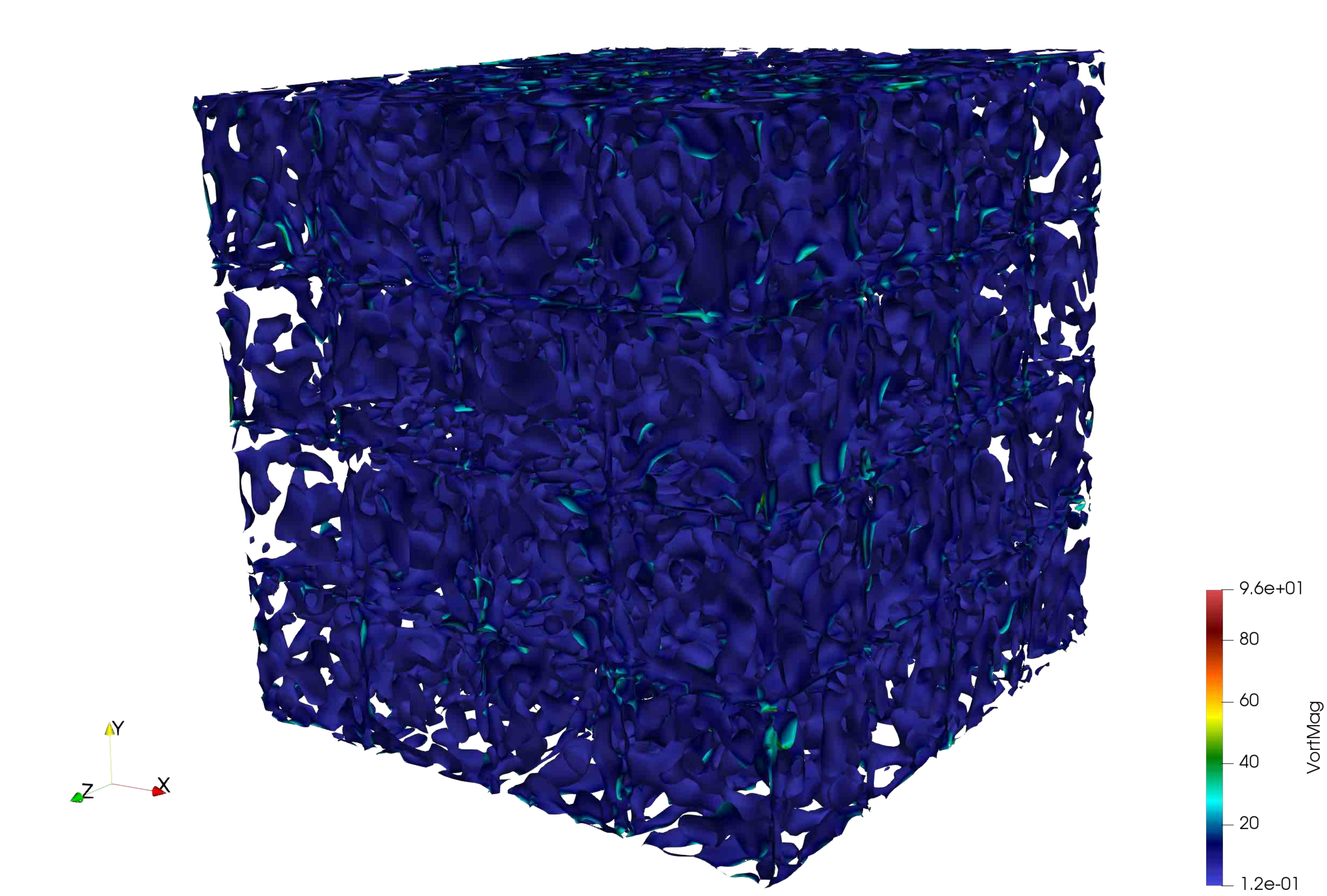}
\end{center}
\caption{Example 3: Transition to decaying homogeneous isotropic 3D-turbulence: 
         5.0 isocontour of Q-criterion, colored with vorticity at 
         $t \in \lbrace 0.0,~6.5,~7.0,~7.5,~10.0,~20.0 \rbrace$}
\label{fig:3Dvortex}
\end{figure}

The snapshots of the solution in Fig.~\ref{fig:3Dvortex} show that until $t \approx 6$, the numerical
method tries to preserve the 2D-behavior of the 2D-lattice flow. This can be seen from the ''vortex tubes'' 
(presented  by the 5.0-isocontour of the so-called Q-criterion, colored with vorticity). Then the vortex 
stretching starts to deform the vortex tubes until $t=7.5$. Later on, i.e. around $t=10$, there starts the 
eddy-breakdown in the inertial range. Here we observe the transition to {\it homogeneous isotropic} 
turbulence.
\begin{figure}[ht]
\begin{center}
\includegraphics[width=3.8cm]{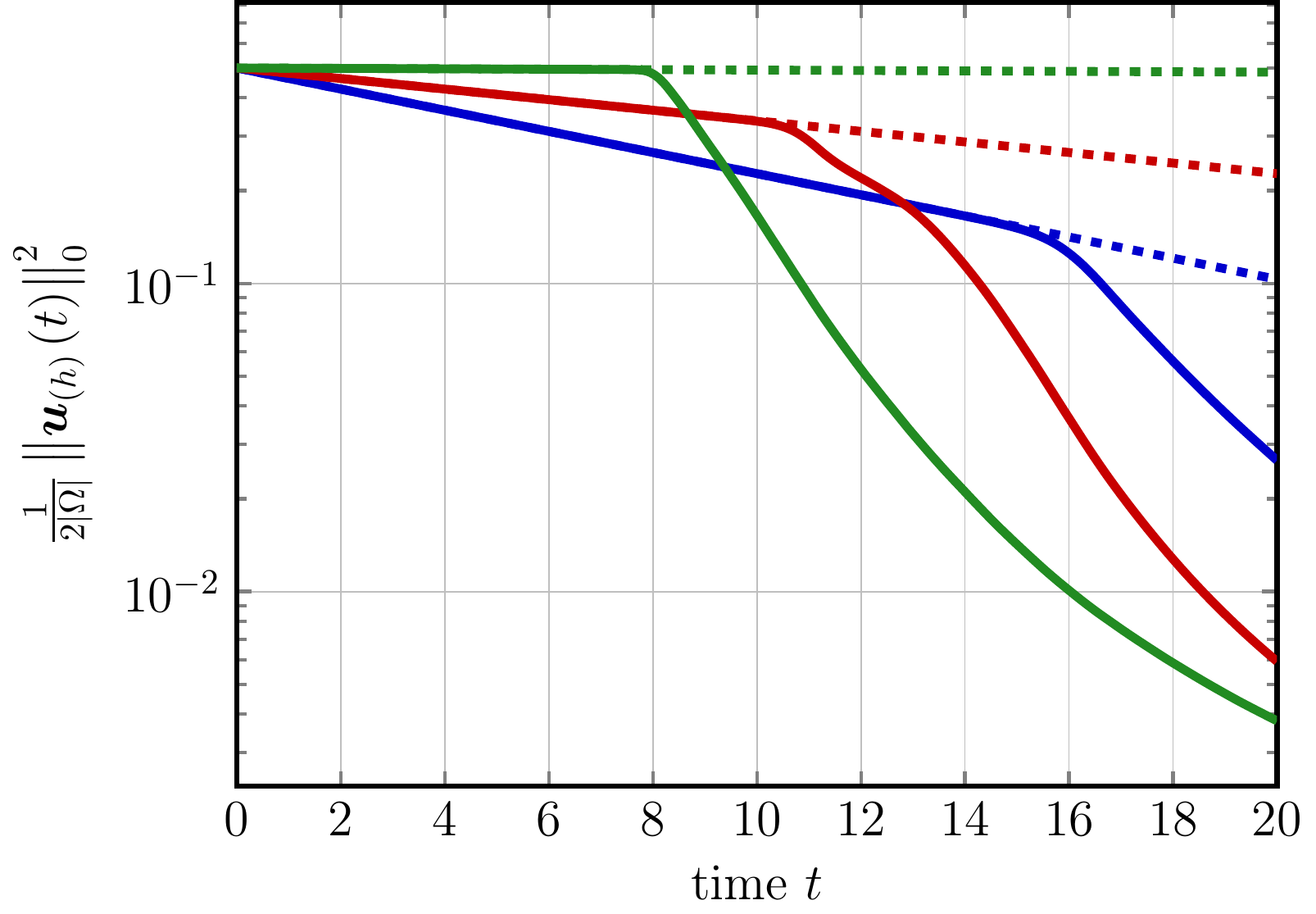} \qquad
\includegraphics[width=3.8cm]{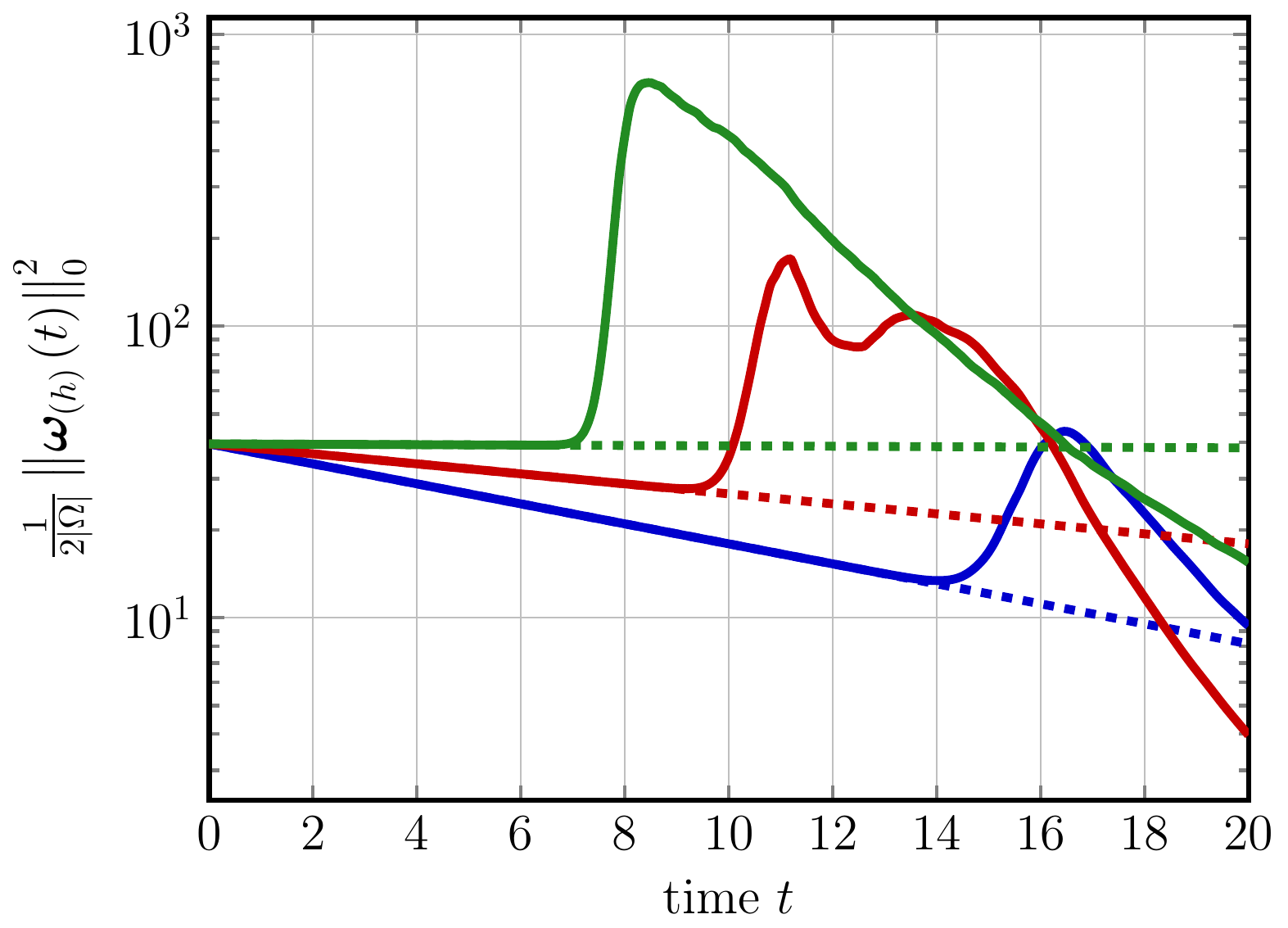} \\
\includegraphics[width=3.8cm]{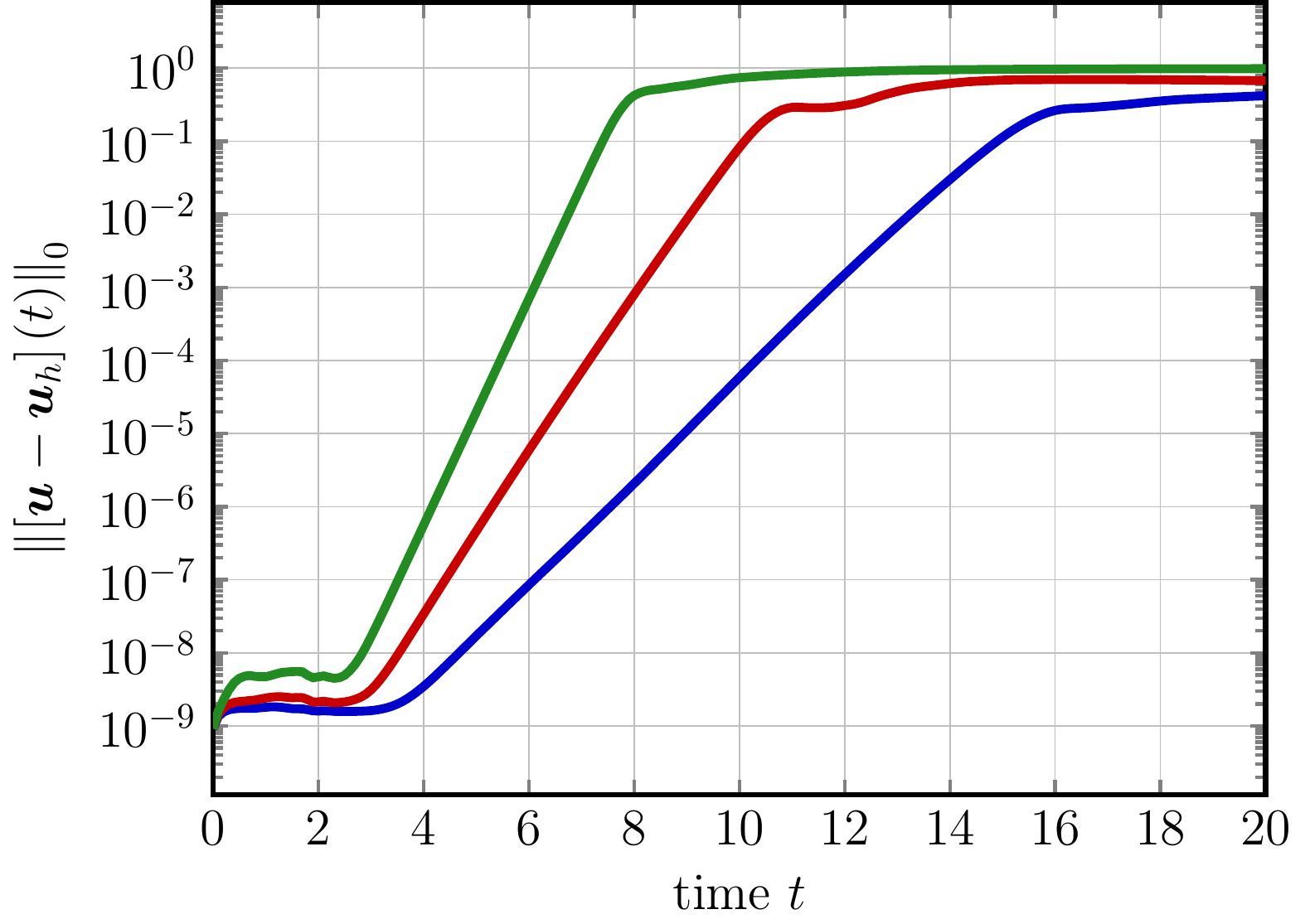} \qquad
\includegraphics[width=3.8cm]{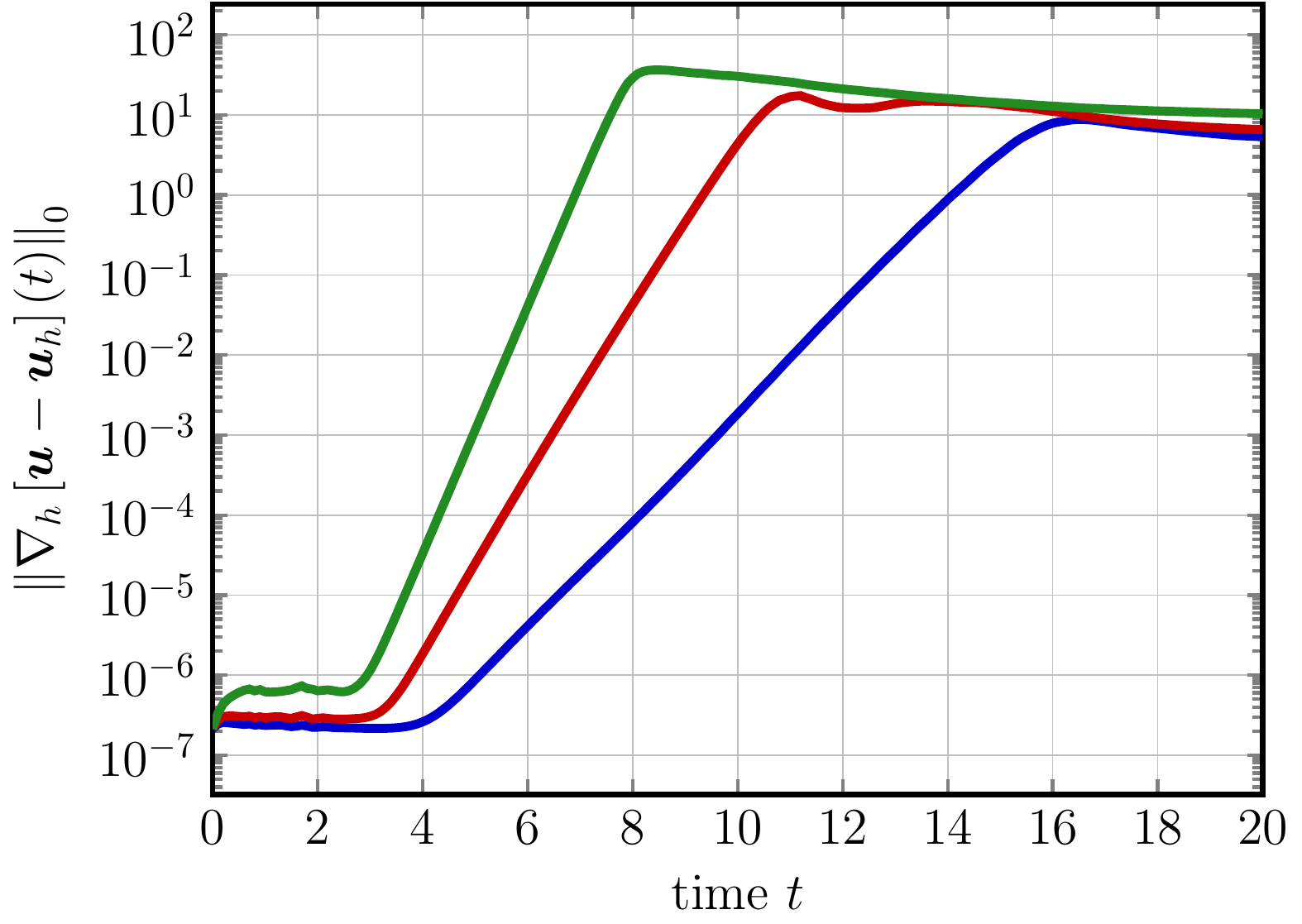} 
\end{center}
\caption{Example 3. First row: $t$-dependent kinetic energy and enstrophy,~
         Second row:  $t$-dependence of errors in $L^2$ and $H^1$, 
         Legend:~ \textcolor{blue}{$\frac{1}{\nu}=2.0 \times 10^3$},
         \textcolor{red}{$\frac{1}{\nu}=4.0 \times 10^3$}, 
         \textcolor{green}{$\frac{1}{\nu}= 1.0 \times 10^5$}  }
\label{fig:3Dvortex-rates}
\end{figure}
Finally, in Fig.~\ref{fig:3Dvortex-rates}, we consider the influence of the Reynolds number for
$\frac{1}{\nu} \in \lbrace 2000, 4000, 100000 \rbrace$.
We apply again the high-order H(div)-dGFEM (here with $k=8, h=\frac{1}{8}$).
In the first row, one observes the strongly decaying kinetic energy and the effect of vortex 
stretching in the (scaled) dissipation rate in time.

The solution is still a Beltrami flow since $(\nabla \times {\bf u}) \times {\bf u} = {\bf 0}$.
Thus a linearization via $p \mapsto P:=p+\frac12 |{\bf u}|^2$ would retain {\it coherent} structures 
as in 2D. This corresponds to the  formal exact solution, see dashed lines. Solid lines correspond to the 
discrete solutions with $k=8$ and $h=\frac18$. The deviation of the discrete solution from the (formal)
exact solution starts earlier for increasing Reynolds number. On the other hand, the deviation can be
shifted to larger times if the FEM-order $k$ is increased and/or an $h$-refinement is performed.

In the second row of Fig.~\ref{fig:3Dvortex-rates}, we consider the $L^2$- and $H^1$-errors for 
${\bf u}-{\bf u}_h$. According to the estimate in Thm.~1, one observes the exponential
behavior of both errors in time. This again indicates that, after a certain time, the discrete solution
deviates from the (formal) exact solution. 
\end{example}
\begin{figure}[ht]
   \centering
     \includegraphics[width=11.0cm]{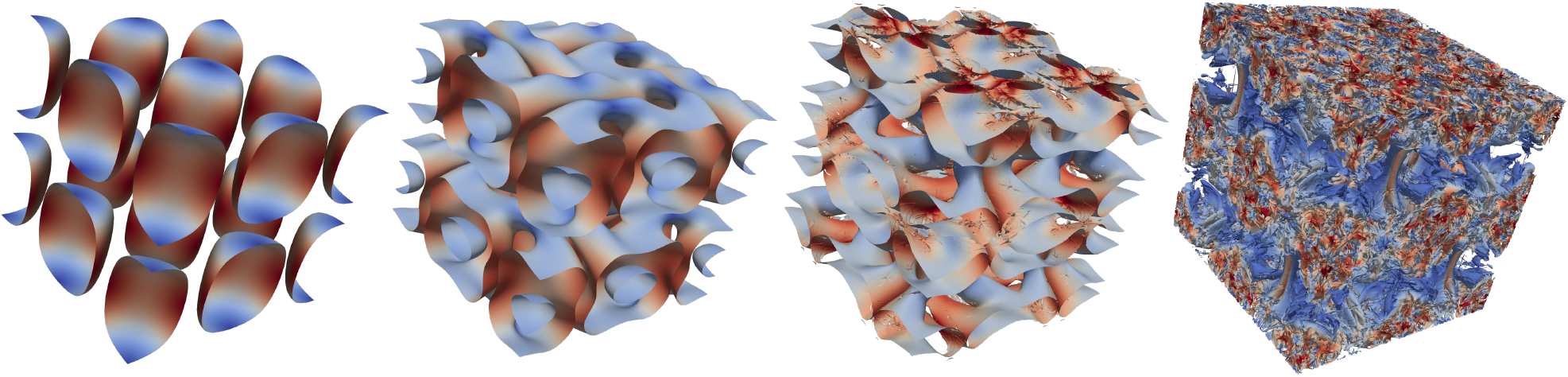}
  \caption{Example 4. Behavior like decaying homogeneous isotropic 3D-turbulence: 
  0.1-isocontour of Q-criterion, colored with velocity at $t \in \lbrace 0.0,~2.0,~4.0,~9.0 \rbrace$}
 \label{fig:Taylor-Green1600}
\end{figure}
\begin{example}~  3D-Taylor-Green vortex at $Re=1600$ 
 \smallskip \\
A typical LES-benchmark is the 3D-Taylor-Green vortex problem at $Re=\frac{UL}{\nu} = 1600$ 
with ${\bf f}={\bf 0}$ and initial condition
$$ {\bf u}_0 (x) = U \Big( \sin \frac{x_1}{L} \cos \frac{x_2}{L} \cos \frac{x_3}{L},
    - \cos \frac{x_1}{L} \sin \frac{x_2}{L} \cos \frac{x_3}{L}, 0 \Big)^t . $$
As in the previous example we observe the breakdown of large eddies into smaller and 
smaller eddies, see Fig.~\ref{fig:Taylor-Green1600}. This indicates that the typical behavior of 
{\it homogeneous isotropic turbulence} develops already for this
relative small Reynolds number $Re=1600$ where we set $U=L=1$. 
\smallskip \\
\begin{figure}[ht]
  \centering
   \includegraphics[width=5.5cm]{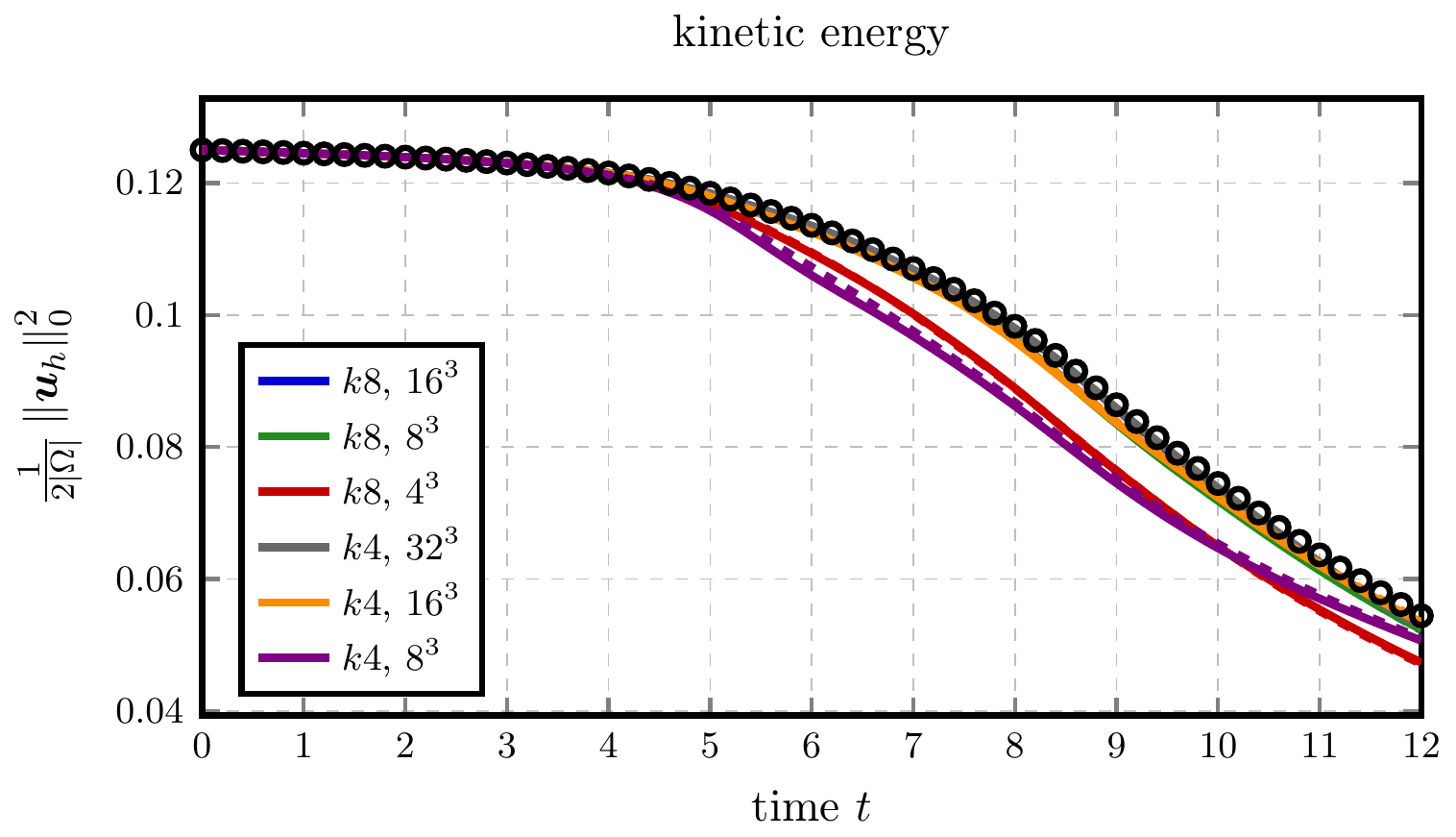} \qquad
    \includegraphics[width=4.5cm]{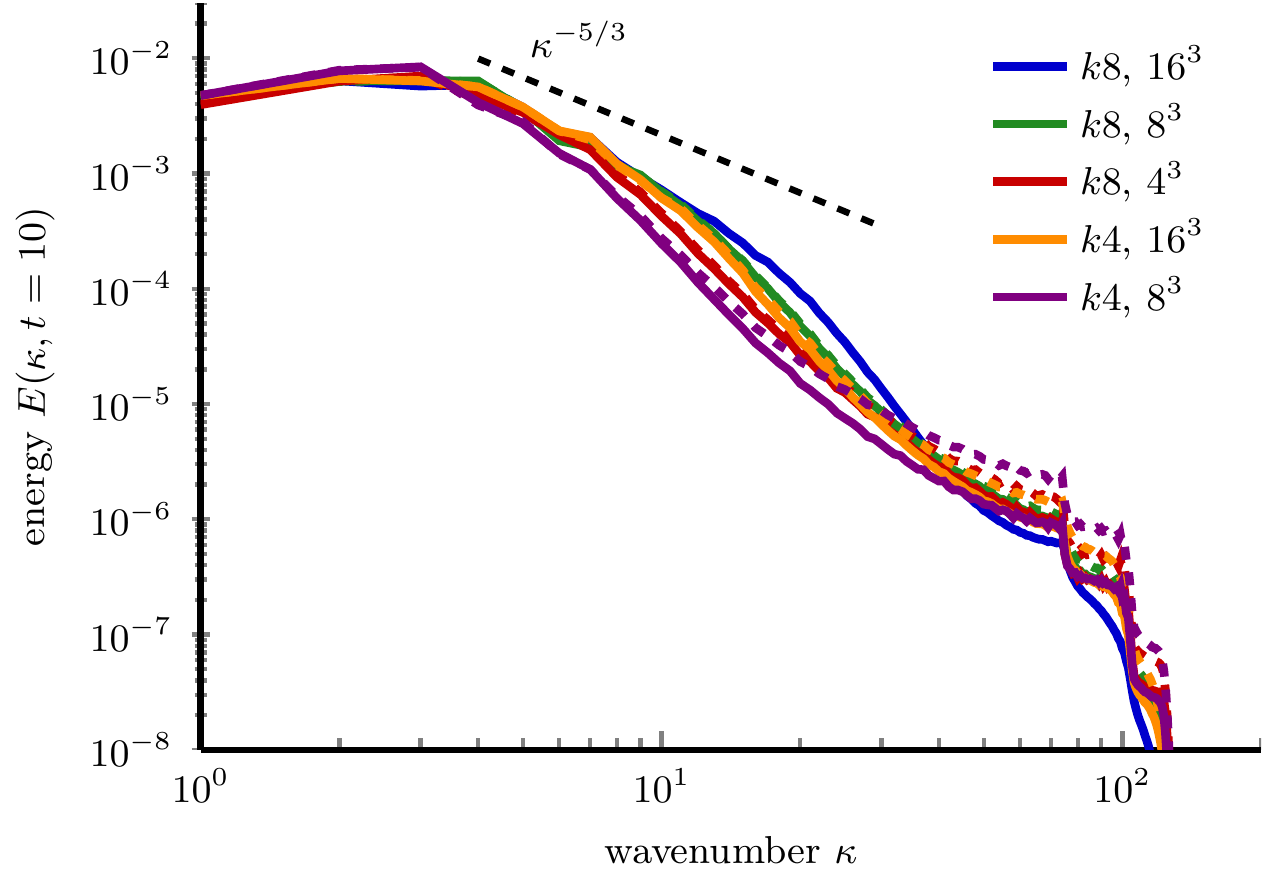} 
  \caption{Example 4. Left: Temporal development of kinetic energy for 3D-Taylor-Green vortex at $Re=1600$ 
           for different values of order $k$ and $N=1/h$;
           Right: Spectrum of kinetic energy at $t=10$ for different values of $k$ and $h=1/N$}
\label{fig:Taylor-Green1600-2}
\end{figure}

Consider now the temporal development of kinetic energy resp. the $L^2$-energy spectrum, see 
Fig.~\ref{fig:Taylor-Green1600-2}.
For ${\bf f}={\bf 0}$, we found in Subsec.~3.1 a weak exponential decay of kinetic energy according
to (\ref{stab1}).
As reference solution serves the solution of a spectral method with $512^3$ grid points (ooo). 
For increasing values of FEM-order $k$ and/or increasing spatial resolution (via refinement of 
$h=1/N$), we observe grid convergence for the kinetic energy, see Fig.~\ref {fig:Taylor-Green1600-2}
(left). 

In Fig.~\ref{fig:Taylor-Green1600-2} (right) we plot the spectra of the kinetic energy at $t=10$ 
for different values of $k$ and $h=1/N$. In particular, no pile-up of the spectra for large 
wave numbers $k$ occurs. The Kolmogorov rate of $E(k) = {\mathcal O}(k^{-5/3})$ is not
reached since $Re=1600$ is too small but will be reached at larger values of $Re$. 
\begin{figure}[ht]
  \centering
  \includegraphics[width=6.0cm]{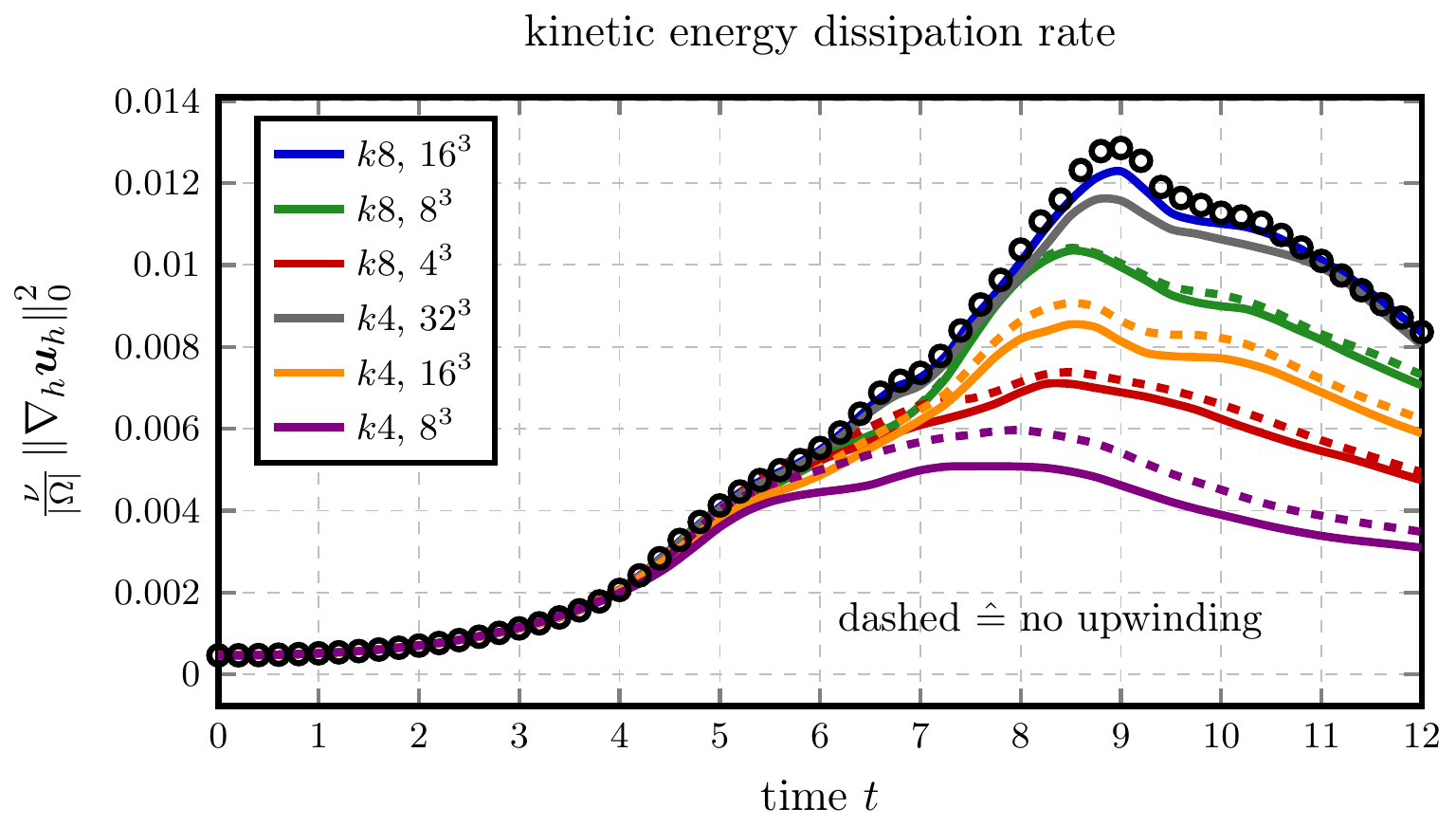} \qquad
   \includegraphics[width=0.4\textwidth]{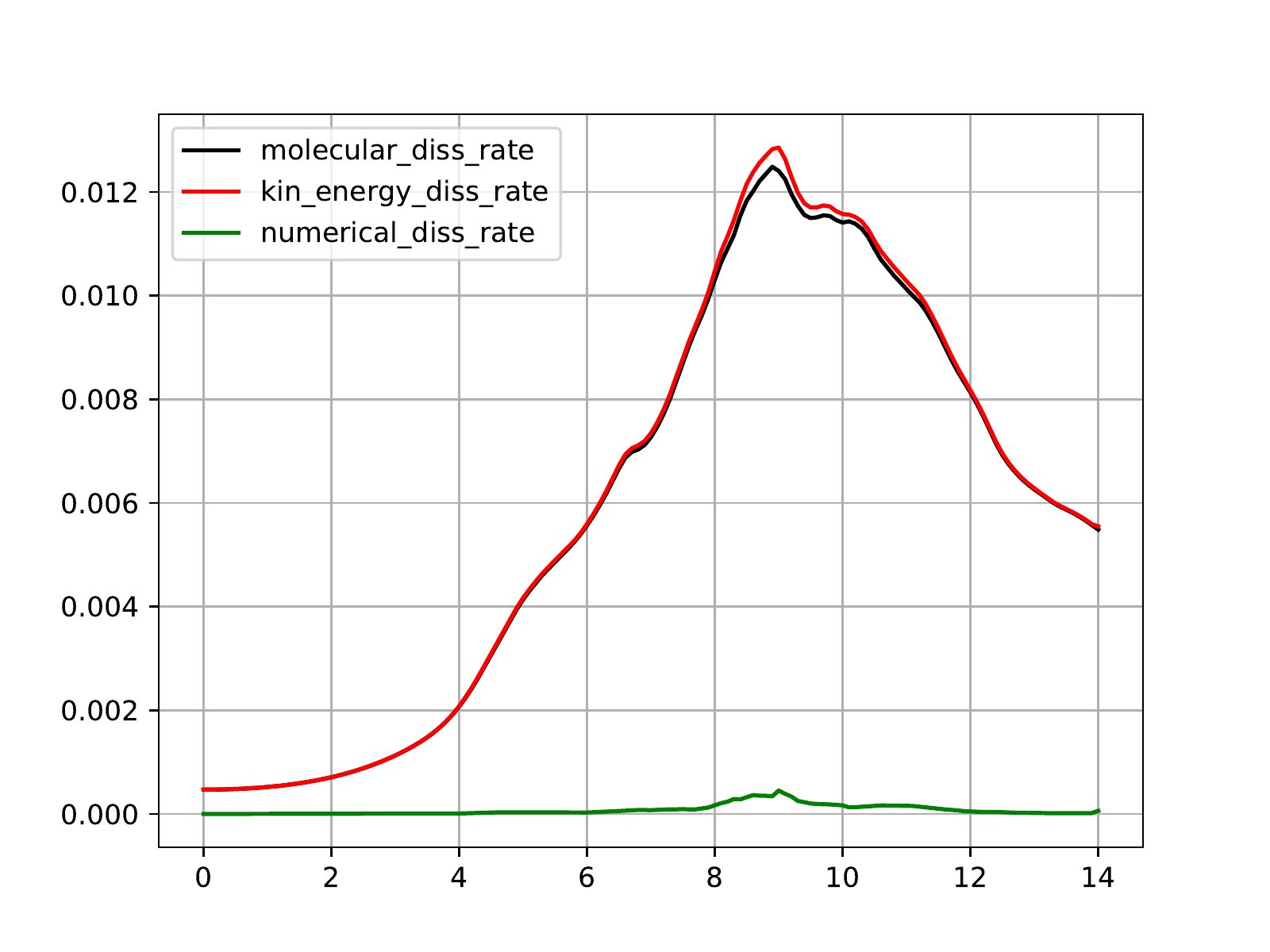} 
   \caption{Example 4. Left: Temporal development of energy dissipation rate for Taylor-Green
            vortex at $Re=1600$ for different values of order $k$ and $N=1/h$;
            Right: Balance of dissipation rates} 
  \label{fig:Taylor-Green1600-4}
\end{figure}

Consider now the temporal development of the kinetic energy dissipation rate for which we obtained 
estimate (\ref{stab2}).
This quantity is much harder to approximate. For increasing values of FEM-order $k$ 
and/or increasing resolution (via refinement of $h=1/N$), we 
observe nearly convergence for the energy dissipation rate. In particular, we find that upwind 
stabilization (see solid lines) decreases the energy dissipation rate on the coarser grids,
see Fig.~\ref{fig:Taylor-Green1600-4} (left). 

Finally, consider the balance of dissipation rates according to
\[
   \frac{d}{dt} \Big( \frac12 \|{\bf u}_h(t)\|_{L^2(\Omega)}^2 \Big) + 
    \nu a_h({\bf u}_h(t),{\bf u}_h(t)) +  |{\bf u}_h(t)|_{{\bf u}_h,\text{upw}}^2   =    0 .
\]
The results are plotted for a relatively fine resolution with order $k=8$ and 
$h=1/N=1/16$. This corresponds to $128^3$ grid points. We observe a very
good agreement between molecular dissipation rate $\nu \| \nabla_h {\bf u}_h\|^2_{L^2(\Omega)}$
and kinetic energy dissipation rate $\frac{d}{dt} \left( \frac12 \|{\bf u}_h(t)\|_{L^2(\Omega)}^2 \right)$,
since the numerical dissipation rate (stemming from SIP penalty and upwinding) reaches not more than 
$3\%$ of the other rates around the peak rate 
of molecular dissipation rate, see Fig.~\ref{fig:Taylor-Green1600-4} (right).

Please note that no explicit turbulence modeling has been applied. The price for such 
results is the H(div)-dGFEM simulation with around $9 \times 10^6$ unknowns.
\end{example}
%

\section{Wall-bounded flow}  \label{sec-3}
%
For wall-bounded turbulent flows, one striking problem is the presence of strong boundary layers, 
e.g. at walls. Another problem is to apply a splitting 
${\bf u}_{(h)} = \langle {\bf u}_{(h)} \rangle +  {\bf u}'_{(h)}$ of the solution
into an averaged velocity with some filter $\langle \cdot \rangle$ and fluctuations.
\begin{example} 3D channel flow \\
%
\begin{figure}[ht] 
   \centering
     \includegraphics[width=0.32\textwidth]{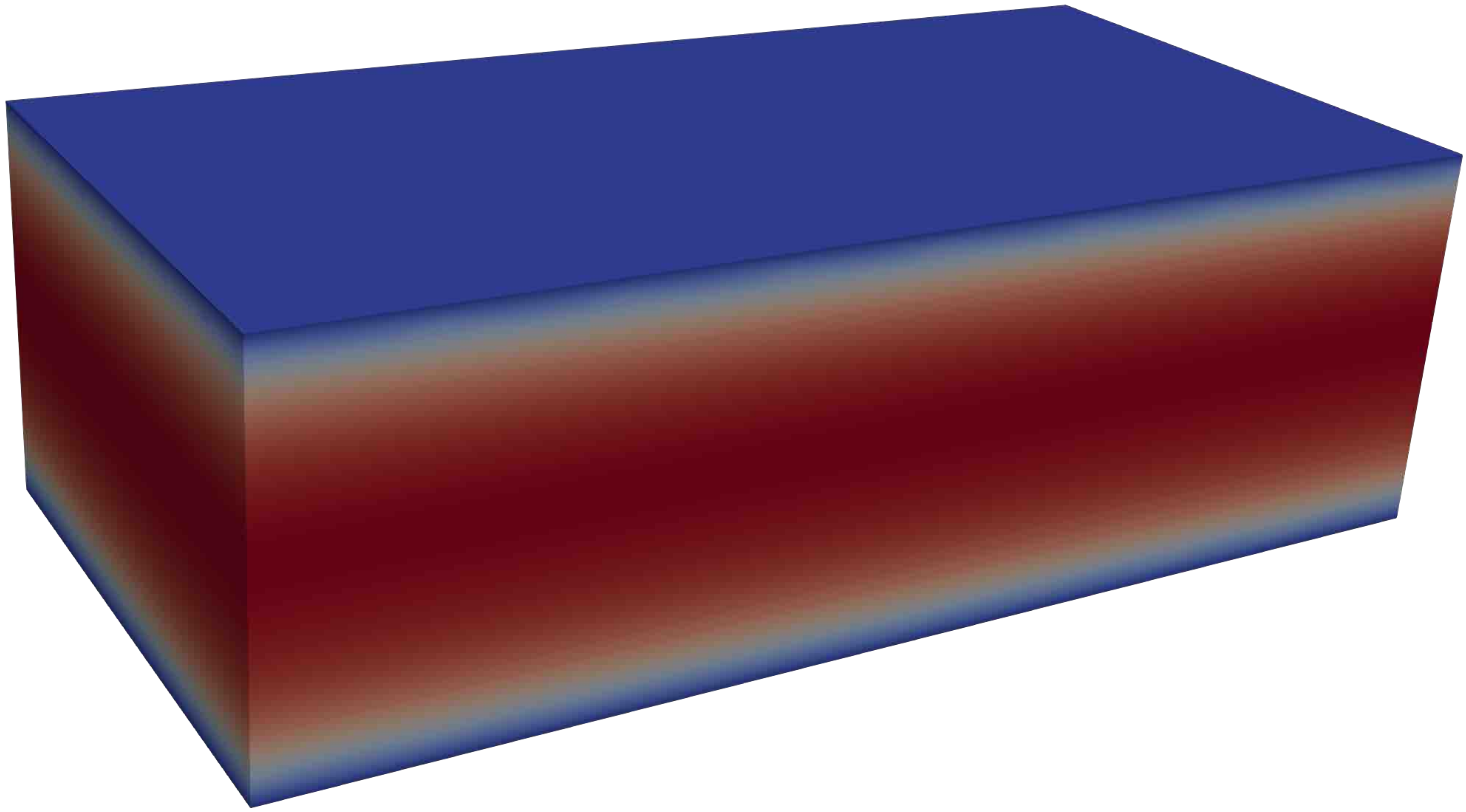} \qquad \quad
     \includegraphics[width=0.32\textwidth]{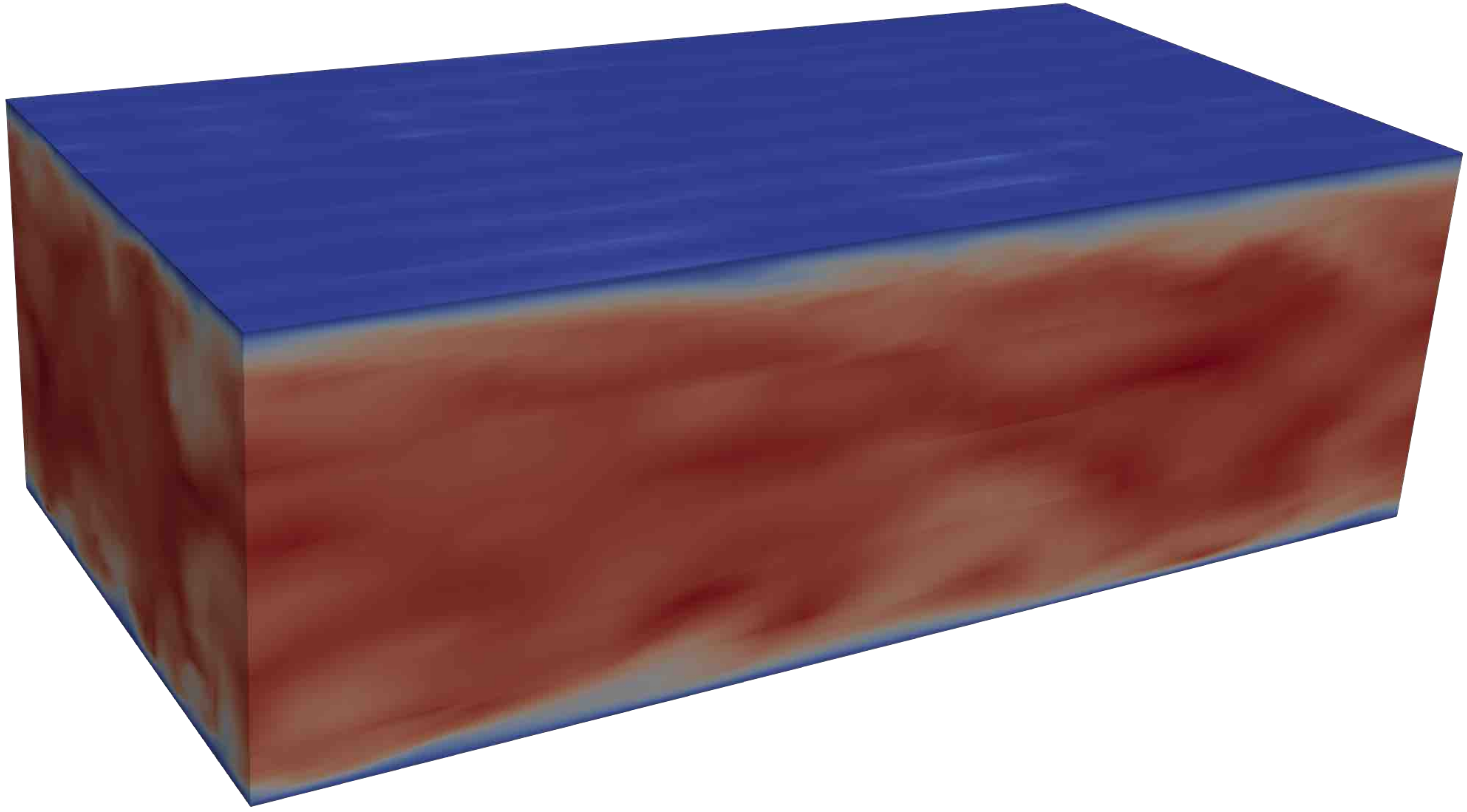} 
  \caption{3D channel flow:~~ Laminar flow (left), \qquad Turbulent $Re_\tau=180$-flow (right)} 
  \label{fig:channel-flow}
\end{figure}
Fig.~\ref{fig:channel-flow} (left) presents a laminar channel flow with a uniquely defined deterministic 
solution. A snapshot of the turbulent channel flow at $Re_\tau =180$ is shown on the right. The
latter is slightly above the transition from laminar to turbulent flow. 
The chaotic solution of turbulent channel flow can be averaged in time and in $x_1$- and $x_3$-directions.
One obtains, after a certain time of averaging, a relatively simple structure of the flow with
$\langle u_1 \rangle = \langle u_1 \rangle (x_2)$.
%
%

Prandtl's boundary layer theory leads to the so-called {\it law of-the-wall}, visible in 
Fig.~\ref{fig:channel180}. 
The mean viscous stress at the wall, the  {\it wall-shear stress}, is 
$ \tau_W = \nu \partial_{x_2} \langle u_1 \rangle|_{x_2=0}.$
An appropriate velocity resp. length-scale in the near-wall region are the friction velocity 
$U_\tau = \sqrt{\tau_W}$ resp. $\eta_\nu = \nu/\sqrt{\tau_W}= \nu/U_\tau$. The
friction-based Reynolds number is defined as $Re_\tau = U_\tau H/\nu$ with channel half 
width $H$. 
The layer can be characterized via the non-dimensional distance from wall in wall units
$x_2^+ = x_2/\eta_\nu = U_\tau x_2/\nu $. 
It is characterized by the viscous wall region $x_2^+ <50$ with dominance of molecular viscosity,
including the steep {\it viscous  sublayer} at the wall with $x_2^+ < 5$, and by the outer 
layer with $x_2^+ > 50$. 

The standard approach to resolve boundary layers is to use a (strongly) anisotropic mesh with 
refinement towards the wall(s). 
Very recent results with a $L^2$-based dGFEM-code by Fehn et al. \cite{Fehn-2018} indicate that 
a strong anisotropic $h$-refinement can be relaxed to a (very) coarse $h$-mesh if higher-order 
FEM are applied.
It turns out that for such (highly) under-resolved turbulent flows a ''medium order'' 
($k=4, \ldots ,8$) is most efficient. Another point is that a purely numerical approach to 
stabilization is applied, i.e. no physical LES or VMS model is used.

\begin{figure}[ht]
 \centering
  \includegraphics[width=0.3\textwidth]{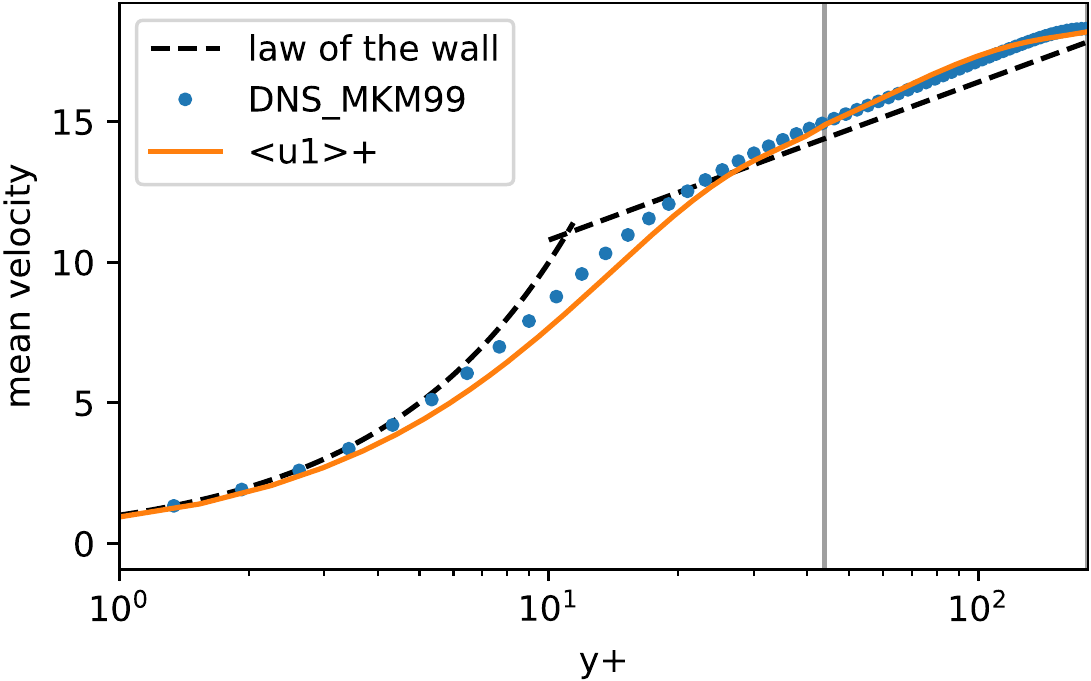} 
   \quad 
  \includegraphics[width=0.3\textwidth]{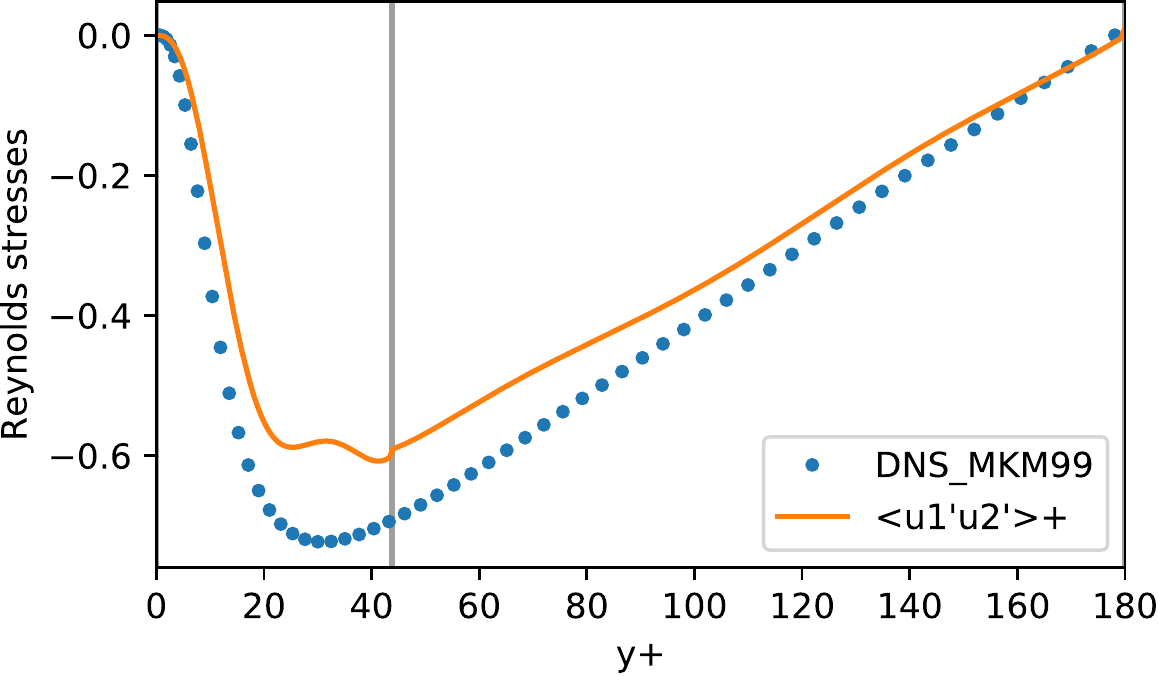}
   \quad 
  \includegraphics[width=0.3\textwidth]{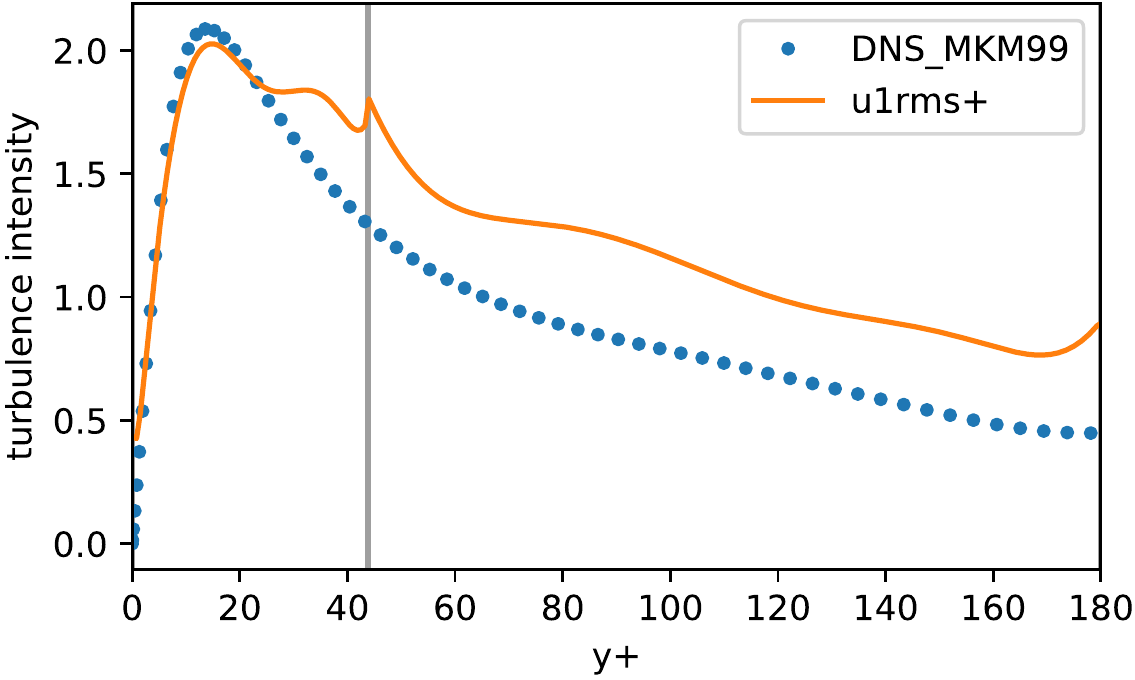} 
 \centering
  \includegraphics[width=0.3\textwidth]{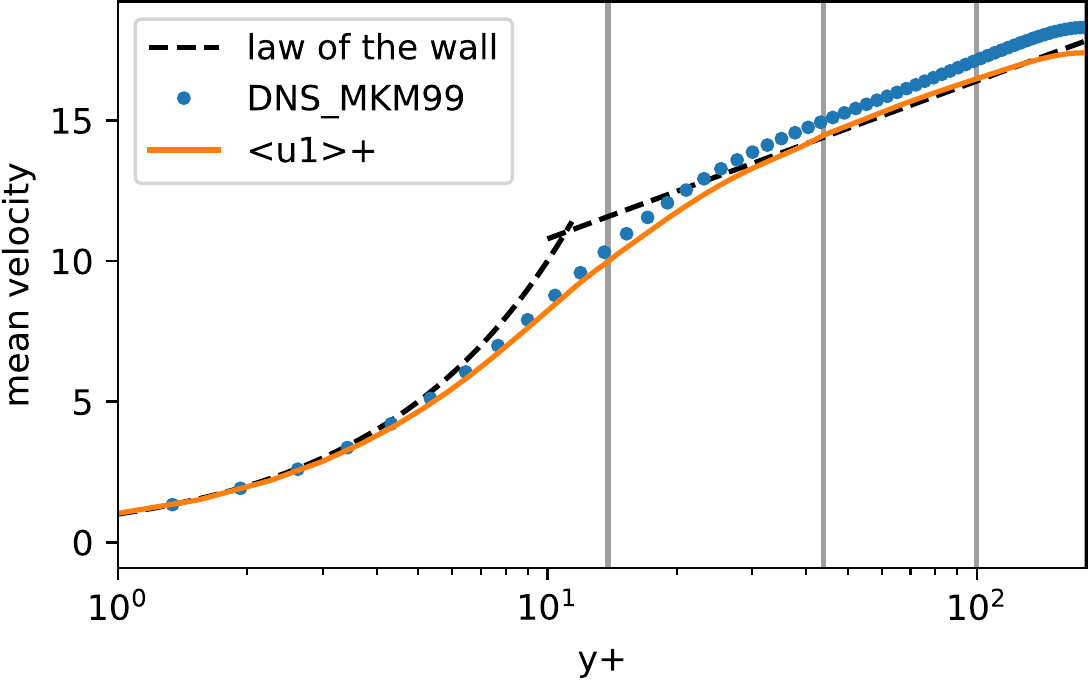} 
    \quad
  \includegraphics[width=0.3\textwidth]{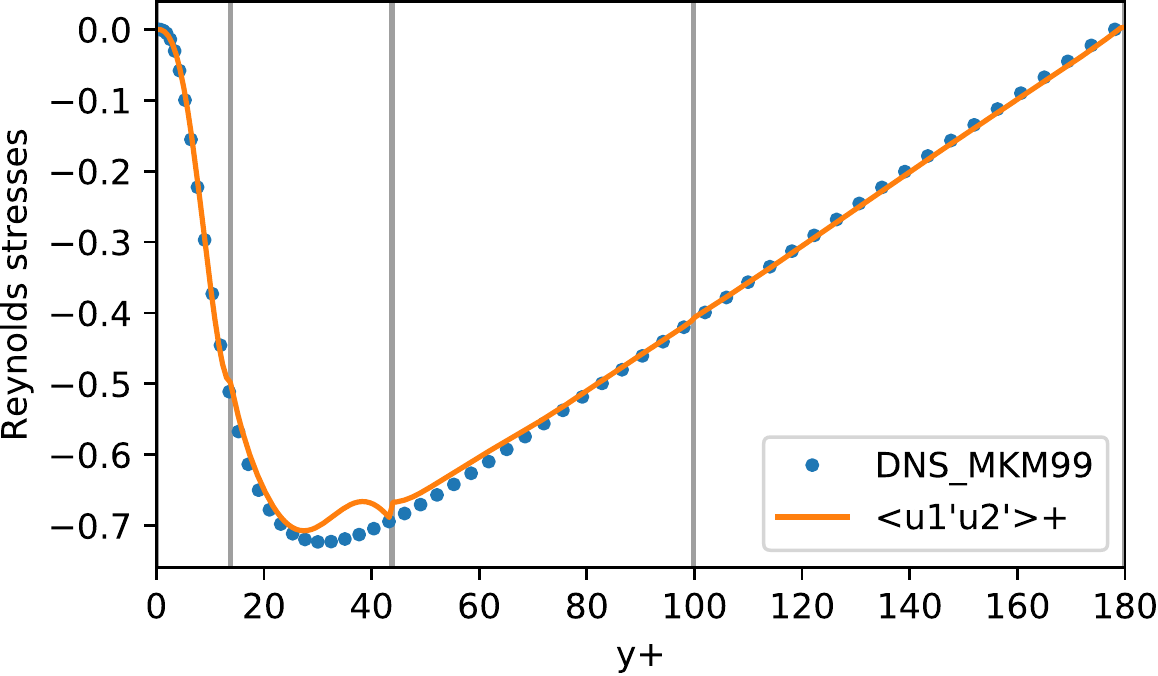} 
   \quad 
  \includegraphics[width=0.3\textwidth]{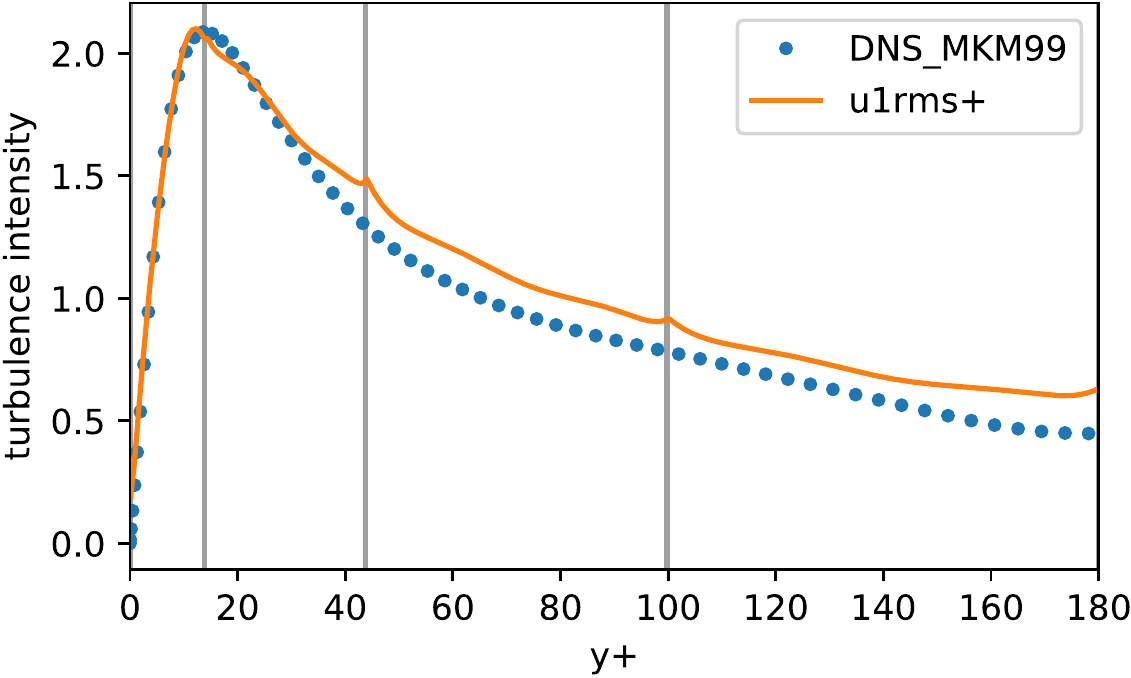} 
 \caption{3D Channel flow at $Re_\tau=180$ with ILES:~ First row: $k=3, N=4$, 
          Second row: $k=3, N=8$; Left: Mean profile $U^+$,~ 
          Middle: Reynolds stress $\langle u'_1u'_2\rangle^+$,~ 
          Right: rms turbulence intensity  $u_{RMS}^+$} 
 \label{fig:channel180}
 \end{figure}

Fig.~\ref{fig:channel180} shows results for the H(div)-dGFEM for the channel flow at 
$Re_\tau = 180$. It turns out that a method of order $k=2$ is not sufficient, but $k=3$
provides good results. Very coarse grids with $N=4$ resp. $N=8$ elements in each 
$x_i$-direction with slightly anisotropic refinement in $x_2$-direction towards the wall
(indicated by vertical lines in Fig.~\ref{fig:channel180}) are used.

No explicit physical LES model is applied. In the ILES approach only numerical 
dissipation (basically from SIP and upwind) is used.  The results for the averaged 
mean profile of  $U^+$, the Reynolds stress component $\langle u'_1u'_2\rangle^+$ and the rms 
turbulence intensity values $u_{RMS}^+$, compared to the DNS-data by Moser et al. \cite{Moser-Kim}, 
are surprisingly good on this very coarse grids with $12^3$ resp. $24^3$ grid points. 
\hfill $\Box$ 
\end{example}
Results in \cite{Fehn-2018} indicate that such approach is also possible for larger 
values of  $Re_\tau$.

\section{Outlook}  \label{sec:4}
%
The following features of H(div)-dGFEM are exploited in the numerical simulation of turbulent flows 
via implicit LES for incompressible Navier-Stokes flows:
 \begin{itemize}
   \item {\it Minimal stabilization:} Numerical dissipation may only result from the SIP term for 
         the diffusive term $a_h$ and upwind term in $c_h$.
   \item {\it Simple form of convective form:} There is no need to modify the convective term 
         $c_h$  since an exactly divergence-free FEM has a clean energy balance {\it a priori}.
   \item {\it Pressure robustness}: H(div)-conforming FEM have the relevant property that
         changing source term ${\bf f}$ to ${\bf f} + \nabla \psi$ changes the solution 
         $({\bf u}_h,p_h)$ to $({\bf u}_h,p_h + \psi)$. 
   \item {\it $Re$-semi-robust error estimates:}  Right-hand-side terms of the error estimate,
         see Thm.~\ref{theorem:converg}, including the Gronwall-term do not explicitly depend on 
         $1/\nu$. 
 \end{itemize}
\begin{figure}[ht]
  \centering
   \includegraphics[width=7.5cm]{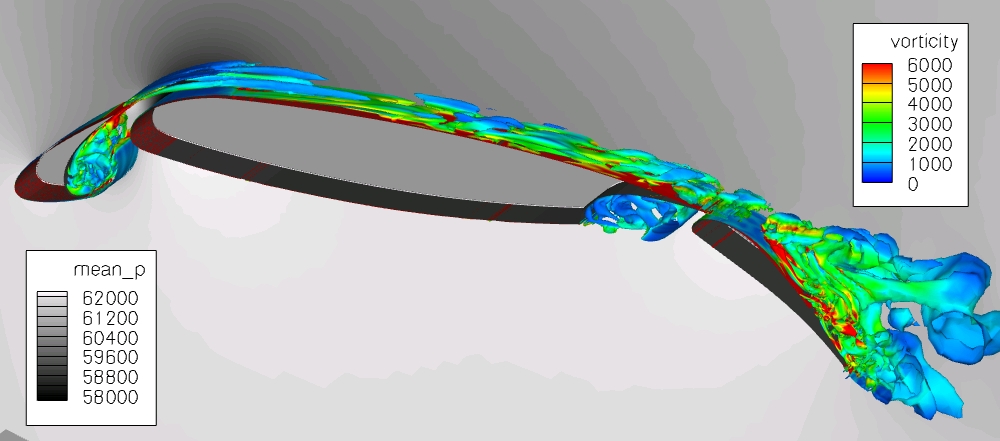}
  \caption{Complex flow around three-element high-lift airfoil}
  \label{fig:3element-profile}
 \end{figure}

We considered an ILES approach to simple turbulent flows with very reasonable results. 
Turbulent flows in practice are clearly much more complex.
A challenge is the flow around a high-lift airfoil, see Fig.~\ref{fig:3element-profile}, 
with complicated interplay of attached laminar and turbulent layers, separation, vortex structures 
etc. For a careful numerical study of such flows see \cite{Reuss-2015}.
A full DNS is still unfeasible. It would be of strong interest to develop new numerical 
concepts for such complex flows which clearly go beyond the limit cases (homogeneous isotropic 
turbulence and turbulent channel flows) under consideration. Nevertheless, the proposed ILES approach 
with high-order and pointwise divergence-free H(div)-dGFEM is a very promising approach. 
Another important point is that the flow in the previous example is governed by the 
{\it compressible} Navier-Stokes model. Many aspects of  incompressible flows can be extended 
to the compressible case, e.g. the approach in boundary layer regions.
%



\end{document}